\documentclass[11pt]{amsart}

\usepackage{amsmath,amsthm,amsfonts,amssymb,bm,graphicx}

\usepackage
{hyperref}
\hypersetup{colorlinks=true,citecolor=blue,linkcolor=blue,urlcolor=blue,
pdfstartview=FitH, pdfauthor=Valentin Blomer ,pdftitle=Applications of the Kuznetsov formula on $GL(3)$}

\theoremstyle{plain}
\newtheorem{theorem}{Theorem}
\newtheorem{lemma}{Lemma}
\newtheorem{cor}[lemma]{Corollary}

\numberwithin{equation}{section}

\theoremstyle{definition}

\renewcommand{\geq}{\geqslant}
\renewcommand{\leq}{\leqslant}

\makeatletter
\DeclareRobustCommand\widecheck[1]{{\mathpalette\@widecheck{#1}}}
\def\@widecheck#1#2{%
    \setbox\z@\hbox{\m@th$#1#2$}%
    \setbox\tw@\hbox{\m@th$#1%
       \widehat{%
          \vrule\@width\z@\@height\ht\z@
          \vrule\@height\z@\@width\wd\z@}$}%
    \dp\tw@-\ht\z@
    \@tempdima\ht\z@ \advance\@tempdima2\ht\tw@ \divide\@tempdima\thr@@
    \setbox\tw@\hbox{%
       \raise\@tempdima\hbox{\scalebox{1}[-1]{\lower\@tempdima\box
\tw@}}}%
    {\ooalign{\box\tw@ \cr \box\z@}}}
\makeatother

\begin{document}

\author{Valentin Blomer}
\author{Jack Buttcane}
\author{P\'eter Maga}
\address{Mathematisches Institut, Bunsenstr. 3-5, 37073 G\"ottingen, Germany} \email{vblomer@math.uni-goettingen.de}
\address{Department of Mathematics, 
244 Mathematics Building, 
University at Buffalo, 
Buffalo, NY 14260-2900 USA} \email{jbuttcane@buffalo.edu}
\address{MTA Alfr\'ed R\'enyi Institute of Mathematics, POB 127, Budapest H-1364, Hungary}
\email{maga.peter@renyi.mta.hu}

\title{Applications of the Kuznetsov formula on ${\rm GL}(3)$ II: the level aspect}

\thanks{The first author was  supported  by the Volkswagen Foundation and a Starting Grant of the European Research Council. The second author was supported by a Starting Grant of the European Research Council. The third author  was supported by a Starting Grant of the European Research Council and  OTKA grant no.\ NK104183.  }


\begin{abstract}  We develop an   explicit Kuznetsov formula on ${\rm GL}(3)$ for congruence subgroups. Applications include a Lindel\"of on average type bound for the sixth moment of ${\rm GL}(3)$ $L$-functions in the level aspect,  an   automorphic large sieve inequality, density results for exceptional eigenvalues and density results for Maa{\ss} forms violating  the Ramanujan conjecture at finite places. 
\end{abstract}


\setcounter{tocdepth}{2}  \maketitle 

 \maketitle

\section{Introduction}

While the toolbox of analytic number theory for classical automorphic forms for congruence subgroups of ${\rm SL}_2(\Bbb{Z})$ is well developed, much less is known in the case of higher rank groups. It is therefore very desirable to extend the collection of available methods by genuine higher rank tools, such as explicit and  for the purpose of analytic number  theory user-friendly spectral summation formulae. In this paper we will introduce a version of the powerful Bruggeman-Kuznetsov formula for congruence subgroups of ${\rm SL}_3(\Bbb{Z})$ and see it in action.  

In the situation of the group ${\rm SL}_2(\Bbb{Z})$ this versatile formula  was first developed independently by Bruggeman \cite{Br} and Kuznetsov \cite{Ku}. Starting with the groundbreaking work of Deshouillers and Iwaniec \cite{DI}, it has become a  very attractive tool, among other things  because it provides a method for studying averages of Kloosterman sums by automorphic techniques, and has shown itself capable of going sometimes beyond the powerful bounds known for individual Kloosterman sums by the Riemann Hypothesis over finite fields. Classical applications include, among others, density results on exceptional eigenvalues, a proof of Selberg's 3/16 theorem, the best known results on the proportion of critical zeros of the Riemann zeta function, and equidistribution of integral points on spheres.

While such a formula exists in great generality, using the theory of general automorphic forms, a main issue for the purpose of analytic number theory is to make the resulting expression analytically useful. This requires a very good understanding of the integral transforms which relate the test functions on both sides of the formulas at all places which turns out to be a problem both in real and $p$-adic analysis. In this paper the focus is on the level aspect, and a good deal of work is devoted to the investigation of the fine properties of ${\rm GL}(3)$ Kloosterman sums with  prime power moduli, but some of the applications also require more precise information on the archimedean test function than those developed  in \cite{Bl}. 

We proceed to describe the new applications that we group into three sections.

\subsection{Moments of $L$-functions} While individual $L$-functions remain   rather elusive objects, statistical information in families $\mathcal{F}$ of $L$-functions is often more easily available. The archetypical result in this direction is a statement on the order of magnitude (asymptotic formulae, upper bounds, or sometimes lower bounds) of some moment of a family of $L$-functions 
$$\sum_{f \in \mathcal{F}} |L(1/2, f)|^k.$$
at the central point. Trace formulae are particularly suitable to evaluate such moments if the family is given by ``spectral properties''.  Here we consider the  ${\rm GL}(3)$ $L$-functions of  large (prime) level $N$ for the congruence subgroup
$\Gamma_0(N)$, the subgroup of matrices in ${\rm SL}_3(\Bbb{Z})$ with bottom row congruent to $(0,  0, \ast)$ modulo $N$ acting on the generalized upper half plane $\Bbb{H}_3$.  This is a subgroup of index $N^{2+o(1)}$ in ${\rm SL}_3(\Bbb{Z})$. For an ${\rm SO(3)}$-invariant subspace of a spherical  cuspidal automorphic representation $\pi \subseteq L_{\text{cusp}}^2(\Gamma_0(N)\backslash \Bbb{H}_3)$  let $\mu_{\pi}$ denote the spectral parameter (Langlands parameter at infinity) of $\pi$, normalized so that it is purely imaginary if the Ramanujan conjecture holds. We fix once and for all a compact set $\Omega \subseteq \mathfrak{a}_{\Bbb{C}}^{\ast}$, the complexified dual of the Lie algebra of the maximal torus in ${\rm PGL}_3(\Bbb{R})$.  If $\Omega$ is not too small, there are roughly $\asymp_{\Omega} N^{2+o(1)}$  such representations $\pi$ with $\mu_{\pi} \in \Omega$. 

It is a fairly straightforward exercise with Kuznetsov formula to prove the following  best-possible (``Lindel\"of-on-average'') bound for the fourth moment:
$$\sum_{\substack{\pi \subseteq L^2(\Gamma_0(N)\backslash \Bbb{H}_3)\\ \mu_{\pi} \in \Omega}} |L(1/2, \pi)|^4 \ll N^{2+\varepsilon}$$
for any $\varepsilon >0$. Here and henceforth in  this paper we will apply the usual $\varepsilon$-convention: the letter $\varepsilon$ denotes an arbitrarily small real number, not necessarily the same on each occurrence. 

The main work of the paper is devoted to bound a \emph{sixth} moment in a best-possible fashion.  This gives us the possibility to highlight the finer details of the Kloosterman side of the Kuznetsov formula and to give a sample argument how to combine this formula with the rest of the machinery of analytic number theory, such as multiple Poisson summation, estimation of multiple character sums and stationary phase type arguments for the archimedean weight functions. 

One of the technical problems is that the Kuznetsov formula -- as any spectral summation formula -- requires a spectrally complete expression. Therefore we have to artificially add the continuous spectrum which unfortunately produces a term of larger of magnitude (the maximal Eisenstein series contribute $N^{5/2}$, see Section \ref{52}). This problem was already faced in the ${\rm GL}(2)$ situation in 
  \cite{DFI} and \cite{BHM};  in \cite{DFI}, a delicate analysis identified a term on the arithmetic side of the Kuznetsov formula that cancelled the continuous spectrum contribution, while in \cite{BHM} this problem was solved by introducing extra zeros in the Mellin transform of the weight function in the approximate functional equation. Both approaches require extremely subtle and precise information on the  archimedean test functions in the Kuznetsov formula that is not easily available in higher rank situations.  Here we deal with this problem by twisting the automorphic forms in question with a fixed character in order to kill the unwanted poles incurred by the continuous spectrum. Our first main theorem is as follows.

\begin{theorem}\label{main} Let $N$ be a large prime and let $p$ be a fixed prime. Let $\chi$ be a primitive character modulo $p$ of order $> 2$ and let  $\Omega \subseteq \mathfrak{a}_{\Bbb{C}}^{\ast}$.  Then
$$\sum_{\substack{\pi \subseteq L^2_{\rm cusp}(\Gamma_0(N)\backslash \Bbb{H}_3)\\ \mu_{\pi} \in \Omega}} |L(1/2, \pi \times \chi)|^6 \ll_{p, \Omega, \varepsilon} N^{2+\varepsilon}$$
for every $\varepsilon > 0$. 
\end{theorem}

\subsection{Spectral mean values and a large sieve} Many applications call for an estimate of Fourier coefficients, averaged over the automorphic spectrum.  For each  $\pi \subseteq L_{\text{cusp}}^2(\Gamma_0(N)\backslash \Bbb{H}_3)$ we choose a newvector $\varpi$ that we normalize such that its Fourier coefficients, defined in \eqref{four} and \eqref{four1} below, satisfy $A_{\varpi}(1, 1) = 1$.  
 The following useful result is the level analogue of 
\cite[Theorem 5]{Bl}:

\begin{theorem}\label{thm2} Let $n, m, N \in \Bbb{N}$, $(mn, N) = 1$, $\Omega \subseteq \mathfrak{a}_{\Bbb{C}}^{\ast}$.  Then we have
$$  \sum_{\substack{\pi \subseteq L_{{\rm cusp}}^2(\Gamma_0(N) \backslash \Bbb{H}_3) \\ \mu_{\pi}\in \Omega}} | A_{\varpi}(n, m)|^2\ll_{\Omega, \varepsilon}  (Nmn)^{\varepsilon}(N^2 +  N^{1/2}nm)$$
for every $\varepsilon> 0$. 
 \end{theorem}

The first term on the right hand side is (up to $\varepsilon$) the number of terms in the sum which dominates the second term provided $nm \ll N^{3/2}$. In particular, in this region the result is best possible and can often be used as a substitute for the Ramanujan conjecture. 

A more refined estimate of this type is the following large sieve inequality for the unramified Hecke eigenvalues $\lambda_{\pi}(n) = A_{\varpi}(n, 1)$, $(n, N) = 1$. 

\begin{theorem}\label{thm3} Let $N \in \Bbb{N}$, $\Omega \subseteq \mathfrak{a}_{\Bbb{C}}^{\ast}$, $X \geq 1$, and let $\alpha(n)$ be a sequence of complex numbers supported on $X \leq n \leq 2X$. Then
$$\sum_{\substack{\pi \subseteq L_{\text{{\rm cusp}}}^2(\Gamma_0(N) \backslash \Bbb{H}_3) \\  \mu_{\pi}\in \Omega}}\Bigl| \sum_{\substack{X \leq n \leq 2X\\ (n, N) = 1}}\lambda_{\pi}(n) \alpha(n)\Bigr|^2\ll_{\Omega, \varepsilon}  (NX)^{\varepsilon} (N^2    +  X^2N^{1/2})\| \alpha \|^2$$
for every $\varepsilon > 0$. 
\end{theorem}

This is   in the spirit of the celebrated large sieve inequalities of \cite{DI}. It should be compared with the case $n=3$ of \cite{Ve} which requires $X  \leq N^{1/4}$ to be optimal whereas our result covers the much larger range $X \leq N^{3/4}$ (cf.\ also \cite{DK} for a different  large sieve inequality).  This shows the advantage of using a powerful tool like the Kuznetsov formula as opposed to the soft methods in \cite{Ve} which on the other hand generalize directly to ${\rm GL}(n)$. 


\subsection{Exceptional eigenvalues and the Ramanujan conjecture} The Ramanujan conjecture is one of the central open problems in the theory of automorphic forms, known only for cohomological forms. In analytic number theory it is often important to control the degree to which the Ramanujan conjecture is violated, and to show that this cannot happen too frequently. The following theorems provide  bounds for the density of forms violating the Ramanujan conjecture at a given place, and we will show in particular that in a quantitative sense almost all Maa{\ss} forms satisfy the Ramanujan conjecture at a given place.  In the eigenvalue aspect this has been investigated in \cite{BBR}. We start with the archimedean place and show the following density result for exceptional Maa{\ss} forms of large level:
\begin{theorem}\label{thm4} Let $N$ be a prime and $\Omega \subseteq \mathfrak{a}_{\Bbb{C}}^{\ast}$.  Then
  $$ \sum_{\substack{\pi \subseteq L_{\text{{\rm cusp}}}^2(\Gamma_0(N) \backslash \Bbb{H}_3) \\ \mu_{\pi} \in \Omega}} N^{4\|\Re  \mu_{\pi}\|} \ll_{\Omega, \varepsilon} N^{2+\varepsilon}$$
  for every $\varepsilon > 0$. 
   \end{theorem}
   Here and in the following, $\| . \|$ denotes the maximum norm. 
The Jacquet-Shalika \cite{JS} bounds imply $\|\Re  \mu_{\pi}\| \leq 1/2$ while the Kim-Sarnak method shows $\|\Re  \mu_{\pi}\| \leq 5/14$. Our result recovers (essentially) the Jacquet-Shalika bounds, but it shows much more: exceptional Maa{\ss} forms occur   less and less frequent, the more the Ramanujan conjecture at infinity is violated. 

A similar result can be obtained for a fixed finite place. Let $\alpha_{\pi}(p )$ 
denote the Satake parameter of a representation  $\pi$ at $p$. The Ramanujan conjecture states that  all three entries of $\alpha_{\pi}( p)$ have absolute value one. 
\begin{theorem}\label{thm5} Let $N \in \Bbb{N}$,  fix a prime $p \nmid N$ and let $\delta > 0$. Let $\Omega \subseteq \mathfrak{a}_{\Bbb{C}}^{\ast}$.  Then there exists $\eta > 0$ (depending on $\delta$ and $p$) such that  
$$\#\left\{ \pi \subseteq L_{\text{{\rm cusp}}}^2(\Gamma_0(N) \backslash \Bbb{H}_3) :  \mu_{\pi} \in \Omega, \, \| \alpha_{\pi}( p)\| 
 \geq 1 + \delta\right\} \ll_{\Omega, \delta, p} (N^2)^{1 - \eta}.$$
\end{theorem}    

\section{The Kuznetsov formula for congruence subgroups of ${\rm SL}_3(\Bbb{Z})$}

In this section we state and prove  the Kuznetsov formula and correct a small error in the statement of the formula in \cite{Bl}. This requires a bit of notational preparation. 
Let $N \in \Bbb{N}$ be the level. We follow the approach in \cite{Bl} and compute the inner product of two Poincar\'e series in two ways.  Let $F : (0, \infty)^2 \rightarrow \Bbb{C}$ be a smooth compactly supported function. Let
\begin{equation}\label{Fast}
  F^{\ast}(y_1, y_2) := F(y_2, y_1). 
\end{equation}
 For two positive integers $m_1, m_2$ and $z  = \left( \begin{smallmatrix} 1 & x_2 & x_3\\ & 1 & x_1\\ & & 1\end{smallmatrix}\right) \left(\begin{smallmatrix} y_1y_2 & & \\ & y_1 & \\ & & 1\end{smallmatrix}\right) \in \Bbb{H}_3$ let $$\mathcal{F}_{m_1, m_2}(z) := e(m_1x_1 + m_2x_2) F(m_1y_1, m_2y_2).$$ Then we consider   the following Poincar\'e series:
\begin{displaymath}
  P_{m_1, m_2}(z) := \sum_{\gamma \in  \Gamma_{\infty} \backslash  \Gamma_0(N)} \mathcal{F}_{m_1, m_2}(\gamma z)
\end{displaymath}
where $\Gamma_{\infty}$ is the subgroup of unipotent upper triangular matrices. The Fourier expansion of these functions features Kloosterman sums and their archimedean analogues, certain special functions given by an integral representation. The three non-trivial terms in the Kuznetsov formula are attached to the elements $w_4 = \left(\begin{smallmatrix} & 1 &\\ & & 1\\ 1 & &\end{smallmatrix}\right)$, $w_5 =  \left(\begin{smallmatrix} &  &1\\ 1& & \\  & 1&\end{smallmatrix}\right)$ and $w_6 =  \left(\begin{smallmatrix} &  &1\\ & -1& \\ 1 & &\end{smallmatrix}\right)$ in the Weyl group. Correspondingly, for $m_1, m_2, n_1, n_2 \in \Bbb{Z} \setminus \{0\}$  we define
\begin{equation}\label{klo45}
	\tilde{S}(m_1,n_1,n_2;D_1,D_2) := \sum_{\substack{C_1 (\text{mod }D_1), C_2 (\text{mod }D_2)\\(C_1,D_1)=(C_2,D_2/D_1)=1}} e\left(n_1\frac{\bar{C_1}C_2}{D_1}+n_2\frac{\bar{C_2}}{D_2/D_1}+m_1\frac{C_1}{D_1}\right),
\end{equation}
for $D_1\mid D_2$   and
\begin{equation}\label{klo6}
\begin{split}
&S^{(N)}(m_1, m_2, n_1, n_2; D_1, D_2)\\
& = \sum_{\substack{B_1, C_1 \, ({\rm mod }\, D_1)\\B_2, C_2 \, ({\rm mod }\,  D_2)\\ D_1C_2 + B_1B_2 + D_2C_1 \equiv 0 \, ({\rm mod }\, D_1D_2)\\ (B_j, C_j, D_j) = 1, N \mid B_1}} e\left(\frac{m_1B_1 + n_1(Y_1 D_2 - Z_1 B_2)}{D_1} + \frac{m_2B_2 + n_2(Y_2 D_1 - Z_2B_1)}{D_2}\right)
\end{split}
\end{equation}
for $N \mid D_1$, $N \mid D_2$, 
where $Y_jB_j + Z_jC_j \equiv 1 \, (\text{mod }D_j)$ for $j = 1, 2$. The latter is almost the same sum as in \cite[Section 4]{BFG} for level 1 except for   the  additional divisibility condition $N \mid B_1$.  Note, however, that   $N\mid B_1B_2$ is automatic (since $N \mid D_1, D_2$), so the additional condition $N \mid B_1$ is relatively minor.

The archimedean functions don't see the additional level and are identical to  the level 1 case.  For $\epsilon \in \{\pm 1\}$ or $\{\pm 1\}^2$,  $F$ as above and $A_1, A_2 > 0$ we define
\begin{equation}\label{test1}
\begin{split}
  \tilde{\mathcal{J}}_{\epsilon;F}(A) & =  A^{-2} \int_0^{\infty}\int_0^{\infty} \int_{-\infty}^{\infty}\int_{-\infty}^{\infty}   e(-\epsilon Ax_1y_1) e\left(y_2 \cdot \frac{x_1x_2}{x_1^2 + 1}\right) e\left(\frac{A}{y_1y_2} \cdot \frac{x_2}{x_1^2 + x_2^2 + 1}\right) \\
& \quad\quad\quad \times F\left(y_2 \cdot \frac{\sqrt{x_1^2+x_2^2 + 1}}{x_1^2 + 1}, \frac{A}{y_1y_2}\cdot \frac{\sqrt{x_1^2 + 1}}{x_1^2+x_2^2 + 1}\right) \overline{F(Ay_1, y_2)} dx_1\, dx_2 \frac{dy_1\, dy_2}{y_1y_2^2},
\end{split}  
\end{equation}
\begin{equation}\label{test2}
\begin{split}
  \mathcal{J}_{\epsilon; F}& (A_1, A_2) =  (A_1A_2)^{-2} \int_0^{\infty}\int_0^{\infty} \int_{-\infty}^{\infty}\int_{-\infty}^{\infty}\int_{-\infty}^{\infty} e\left(-\epsilon_1 A_1x_1y_1 - \epsilon_2 A_2x_2y_2\right)\\
 & \times  e\left(-\frac{A_2}{y_2} \cdot \frac{x_1x_3+x_2}{x_3^2 + x_2^2 + 1}\right) e\left(-\frac{A_1}{y_1} \cdot \frac{x_2(x_1x_2-x_3) + x_1} {(x_1x_2-x_3)^2 + x_1^2 + 1}\right)\overline{F(A_1y_1, A_2 y_2)} \\
& \times F\left(\frac{A_2}{y_2} \cdot \frac{\sqrt{(x_1x_2-x_3)^2 + x_1^2 + 1}}{x_3^2 + x_2^2 + 1}, \frac{A_1}{y_1} \cdot \frac{ \sqrt{x_3^2 + x_2^2 + 1}} {(x_1x_2-x_3)^2 + x_1^2 + 1}\right)  dx_1\, dx_2\, dx_3 \frac{dy_1\, dy_2}{y_1y_2}.
\end{split}  
\end{equation}
Next we define for $\mu \in \mathfrak{a}_{\Bbb{C}}^{\ast}$ and $y_1, y_2 > 0$ 
the (slightly renormalized) Whittaker function as in \cite[(2.15)]{Bl} by its double Mellin transform
\begin{equation}\label{whit}
\begin{split}
  \tilde{W}_{\mu}(y_1, y_2) & = \frac{y_1y_2 \pi^{\frac{3}{2}  } }{     |\Gamma(\frac{1}{2} (1 + i \Im (\mu_1 + 2\mu_2)))\Gamma(\frac{1}{2} (1 + i \Im (\mu_1 - \mu_2)))\Gamma(\frac{1}{2} (1 + i \Im (2\mu_1 + \mu_2)))|}\\
  & \times \frac{1}{(2\pi i)^2} \int_{(1)} \int_{(1)} \frac{\prod_{j=1}^3 \Gamma(\frac{1}{2}(s_1 + \mu_j)) \prod_{j=1}^3 \Gamma(\frac{1}{2}(s_2 - \mu_j)) }{4\pi^{s_1+s_2} \Gamma(\frac{1}{2}(s_1+s_2))} y_1^{-s_1} y_2^{-s_2} ds_1\, ds_2. 
  \end{split}
\end{equation}
For a (not necessarily cuspidal) automorphic form    $\varpi $ of level $N$ and spectral parameter $\mu$ we define the Fourier coefficient $\tilde{A}_{\varpi}(m_1, m_2)$ ($m_1, m_2 \not= 0$) by 
\begin{equation}\label{four}
  \int_0^1\int_0^1 \int_0^1 \varpi(z) e(-m_1x_1-m_2x_2) dx_1 \, dx_2 \, dx_3 = \frac{\tilde{A}_{\varpi}(m_1, m_2)}{|m_1m_2|} \tilde{W}_{\mu}(|m_1|y_1, |m_2|y_2). 
\end{equation}
To ease notation, we will denote by $\{\varpi\}$ an orthonormal basis of automorphic forms of level $N$, cuspidal or Eisenstein series, containing all cuspidal  newvectors, and we denote by $\int_{(N)} d\varpi$ a combined sum/integral over the complete spectrum of level $N$. The relevant spectral decomposition is a special case of Langlands' general theory, see e.g.\ \cite{Ar} for a convenient summary in adelic language. By Hecke theory, we can and will assume that all $\varpi$ are eigenfunctions of the Hecke algebra coprime to $N$. Since 
$$\Gamma_0(N)\text{diag}(m_0m_1m_2, m_0m_1, m_0) \Gamma_0(N) = \Gamma_0(1)\text{diag}(m_0m_1m_2, m_0m_1, m_0) \Gamma_0(1)$$ 
for  $(m_0m_1 m_2, N) = 1$, this is just the unramified Hecke algebra that satisfies the  usual ${\rm GL}(3)$ Hecke relations as in \cite[Theorem 6.4.11]{Go}. The proof of \cite[Theorem 6.4.11]{Go} also shows that if $\tilde{A}_{\varpi}(1, 1) = 0$, then $\tilde{A}_{\varpi}(m_1, m_2) = 0$ whenever $(m_1m_2, N) =1$. If $\tilde{A}_{\varpi}(1, 1)\not= 0$, which  is the case in particular for newvectors $\varpi$, we write 
\begin{equation}\label{four1}
A_{\varpi}(m_1, m_2) = \tilde{A}_{\varpi}(m_1, m_2)/\tilde{A}_{\varpi}(1, 1),
\end{equation}
 in which case the normalized Fourier coefficients $A_{\varpi}(m_1, m_2) $ satisfy the multiplicativity relations of \cite[Theorem 6.4.11]{Go}. If $\tilde{A}_{\varpi}(1, 1)= 0$, we simply write $A_{\varpi}(m_1, m_2) = \tilde{A}_{\varpi}(m_1, m_2)$ and remark already at this place that for such $\varpi$ only vanishing Fourier coefficients will come up in our analysis (which, in a trivial way, satisfy the Hecke relations), so that the normalization is irrelevant. For notational consistency we write $N(\varpi) = \tilde{A}_{\varpi}(1, 1)^{2}$ if $\tilde{A}_{\varpi}(1, 1) \not=0$ and $N(\varpi) = 1$ otherwise.

Rankin-Selberg theory shows   (see e.g.\ \cite[Lemma 1]{Bl}) that for a cuspidal newform $\varpi \in \pi$ one has 
\begin{displaymath}
\mathcal{N}(\varpi)   \asymp [{\rm SL}_3( \Bbb{Z}) : \Gamma_0(N)] \cdot \underset{s=1}{\text{res}} \sum_{m_1, m_2} \frac{|A_{\varpi}(m_1, m_2)|^2}{m_1^{2s} m_2^s},
\end{displaymath}
and  it  follows from  \cite[Theorem 2]{Li} that 
\begin{equation}\label{res}
  \mathcal{N}(\varpi) \ll N^2(N(1 + |\mu_{\pi}|))^{\varepsilon}. 
  \end{equation}
 We define an inner product on  $(0, \infty)^2$ by
\begin{displaymath}
  \langle f, g \rangle := \int_0^{\infty} \int_0^{\infty}  f(y_1, y_2) \overline{g(y_1, y_2)} \frac{dy_1 dy_2}{(y_1y_2)^3}. 
\end{displaymath}
With this notation we are ready to state our version of the Kuznetsov formula. 
\begin{theorem}\label{kuz} Let $F$ be a compactly supported test function with $F^{\ast}$ as in \eqref{Fast}. Let $N, n_1, n_2, m_1, m_2 \in \Bbb{N}$.  Then
\begin{equation}\label{kuzform}
\int_{(N)} \frac{\overline{A_{\varpi}(n_1, n_2)} A_{\varpi}(m_1, m_2)}{\mathcal{N}(\varpi)} |\langle \tilde{W}_{\mu_{\pi}}, F\rangle|^2  d\varpi= \Delta + \Sigma_4 + \Sigma_5 + \Sigma_6
\end{equation}
where
\begin{equation}\label{terms}
\begin{split}
  \Delta& = \delta_{n_1, m_1} \delta_{n_2, m_2}   \| F \|^2,\\
   \Sigma_{4}& = \sum_{\epsilon  = \pm 1} \sum_{\substack{N D_2 \mid D_1\\  n_2 D_1= m_1 D_2^2}}\frac{ \tilde{S}(\epsilon m_2, n_2, n_1, D_2, D_1)}{D_1D_2}  \tilde{\mathcal{J}}_{\epsilon; F^{\ast}}\left(\sqrt{\frac{n_1 n_2 m_2}{D_1D_2}} \right),\\
 \Sigma_{5} &= \sum_{\epsilon  = \pm 1} \sum_{\substack{ N \mid  D_1 \mid D_2\\ n_1 D_2 = m_2 D_1^2}} \frac{ \tilde{S}(\epsilon m_1, n_1, n_2, D_1, D_2) }{D_1D_2}\tilde{\mathcal{J}}_{\epsilon; F}\left( \sqrt{\frac{n_1 n_2 m_1}{D_1D_2}}  \right), \\
   \Sigma_6 &= \sum_{\epsilon_1, \epsilon_2 = \pm 1} \sum_{N \mid D_1, N \mid D_2  } \frac{S^{(N)}(\epsilon_2 m_2, \epsilon_1 m_1, n_1, n_2, D_1, D_2)}{D_1D_2}\mathcal{J}_{\epsilon;F}\left( \frac{\sqrt{n_2 m_1 D_1}}{D_2 }, \frac{\sqrt{n_1 m_2 D_2}}{D_1}\right).
\end{split}
\end{equation} 
\end{theorem}

\textbf{Remarks:} (1) In \cite{Bl}, the first two entries in the long Weyl element Kloosterman sum are mistakenly interchanged, cf.\ \cite[p.\ 64]{BFG}. 



(2) The Fourier coefficients of Eisenstein series for $\Gamma_0(N) \subseteq {\rm SL}_3(\Bbb{Z})$ for \emph{all} indices are computed  in detail in \cite{Ba}. 


(3) Note that there is a small asymmetry in the definition of $\Sigma_4$ and $\Sigma_5$. If $(n_1m_1, N) = 1$, then the summation condition in $\Sigma_4$ is equivalent to $D_1 = N^2d_1d_2$, $D_2 = Nd_2$, $n_2d_1 = m_1d_2$, while the summation condition in $\Sigma_5$ is equivalent to $D_1 = Nd_1 $, $D_2 = N^2d_1d_2$, $n_1d_2 = m_2d_1$, so complete symmetry between $\Sigma_4$ and $\Sigma_5$ is restored. \\

\textbf{Proof.} This is exactly as in \cite[Proposition 4]{Bl} by computing the inner product of two level $N$ Poincar\'e series by unfolding and by spectral decomposition and then comparing both expressions. 
 We only have to verify that the definition of our Kloosterman sums agrees with the Fourier expansion of the level $N$ Poincar\'e series.
The exponential sums appearing in the latter are most easily, but abstractly, defined in terms of the Bruhat decomposition, so the procedure is to enumerate the terms in the sum using the Pl\"ucker coordinates, determine the summand as a function of the Pl\"ucker coordinates by writing out the Bruhat decomposition of each term, and then verify that the summand only depends on the residue classes of the Pl\"ucker coordinates.

Let $U(R)$ be the group of upper-triangular matrices with ones on the diagonal and entries in the ring $R$, $W$ the Weyl group and $V$ the diagonal orthogonal matrices of ${\rm SL}_3(\Bbb{Z})$.
We also need the decomposition of $U(R)$ by $w\in W$, so set $$U_w(R) = (w^{-1} U(R) w) \cap U(R), \quad \bar{U}_w(R) = (w^{-1} U(R)^t w) \cap U(R).$$
Define characters of $U(\Bbb{R})$ by $$\psi_{n_1,n_2}\left(\begin{matrix}1&x_2&*\\&1&x_1\\&&1 \end{matrix}\right)=e\left(n_1 x_1 +n_2 x_2\right)$$ where we assume $n_1,n_2\in\Bbb{Z}$.
Then the Bruhat decomposition of some $\gamma\in{\rm SL}_3(\Bbb{Z})$ takes the form $\gamma=bcvwb'$ with $w\in W$, $v\in V$, $b,b'\in U(\Bbb{R})$ and $c = \text{diag}(1/c_2, c_2/c_1, c_1)$ 
for some $c_1,c_2\in\Bbb{N}$.  The Bruhat decomposition is only defined up to an element of $U_w(\Bbb{R})$.

Now let $w \in W$, $n_1, n_2, m_1, m_2 \in \Bbb{Z}$ and $c = \text{diag}(1/c_2, c_2/c_1, c_1)$ as before. If the \emph{compatibility condition} $$\psi_{n_1,n_2}((cw) u (cw)^{-1}) \psi_{m_1,m_2}(u^{-1}) = 1 \quad \text{for all }   u\in U_w(\Bbb{R})$$ holds, we define  the Kloosterman sums
\begin{equation}\label{kloabs}
S_w(\psi_{n_1,n_2},\psi_{m_1,m_2}; c) = \sum_{\gamma = bcwb' \in U(\Bbb{Z})\backslash\Gamma_0(N)/V\bar{U}_w(\Bbb{Z})} \psi_{n_1,n_2}(b) \psi_{m_1,m_2}(b'). 
\end{equation}
The sum is over representatives $\gamma$ in the quotient space having the prescribed components $c$ and $w$ in their Bruhat decomposition, which is well-defined by the compatibility condition. The quotient by $V$ simply allows us to restrict to positive moduli $c_1$ and $c_2$ by conjugating the $v$ matrix, which contains the signs of the moduli, to the right. If the compatibility relation fails, we simply define $S_w(\psi_{n_1,n_2},\psi_{m_1,m_2}; c) =0$.   

By a computation of Friedberg \cite[pp.\ 173-174]{FrPoincare}, only sums satisfying the compatibility condition occur in the Fourier expansion of a Poincar\'e series. In particular, for $n_1 n_2 m_1 m_2 \ne 0$, only the $I$, $w_4$, $w_5$, and $w_6$ Weyl elements contribute, since otherwise the compatibility relation is never satisfied. 

We now wish to show that the concrete expressions for the Kloosterman sums given in \eqref{klo45} and \eqref{klo6} match the abstract definition \eqref{kloabs}.

We may parameterize representatives of $U(\Bbb{Z})\backslash\Gamma_0(N)$ by the Pl\"ucker coordinates $A_1,B_1,C_1$ and $A_2,B_2,C_2$ satisfying
\begin{align}
\label{eq:PlukerConds}
\begin{array}{l}
	(A_1,B_1,C_1)=(A_2,B_2,C_2)=1, \\
	A_1 C_2+B_1B_2+C_1 A_2=0, \\
	N\mid A_1,\,  N\mid B_1.
\end{array}
\end{align}
For a matrix $\gamma = \left(\begin{smallmatrix}g&h&i\\d&e&f\\a&b&c\end{smallmatrix}\right)\in \Gamma_0(N)$, these are computed by
\[ A_1=a, \quad B_1=b, \quad C_1=c, \quad A_2=bd-ae, \quad B_2=af-cd, \quad C_2=ce-bf. \]
Our computation now essentially follows \cite{BFG}, but we must keep track of the level condition $N\mid A_1, N\mid B_1$.
The auxiliary parameters $Z_2=g$, $Y_2=h$, $X_2=i$, $Z_1=ge-dh$, $Y_1=di-gf$, and $X_1=fh-ei$ are solutions to the equations
\begin{align}
\label{eq:AuxPluker}
\begin{array}{l}
	Z_2 C_2+Y_2 B_2+X_2 A_2=1, \\
	Z_1 C_1+Y_1 B_1+X_1 A_1=1.
\end{array}
\end{align}
These equations do not completely determine the auxiliary parameters, but we will only require that the auxiliary parameters are some solution, as the final expression for the Kloosterman sum will be independent of the choice. 
The right-translation action of $x \in U(\Bbb{Z})$ on the Pl\"ucker coordinates gives the new values
\begin{align*}
\begin{array}{l}
(A_1, B_1,C_1) \mapsto (A_1, B_1+x_2 A_1,C_1+x_1 B_1+x_3 A_1), \\
(A_2, B_2, C_2) \mapsto (A_2, B_2 - x_1 A_2, C_2 - x_2 B_2 + (x_1 x_2-x_3) A_2).
\end{array}
\end{align*}

Now the Bruhat decomposition for elements of the long element Weyl cell in $\Gamma_0(N)$ may be written as
\[ \gamma = \left(\begin{matrix}1&\frac{Z_2 B_1 - Y_2 A_1}{A_2}&\frac{Z_2}{A_1} \\&1&\frac{Y_1 A_2-Z_1 B_2}{A_1}\\&&1\end{matrix}\right) \left(\begin{matrix}\frac{1}{A_2}\\&\frac{A_2}{A_1}\\&&A_1\end{matrix}\right) w_6 \left(\begin{matrix}1&\frac{B_1}{A_1}&\frac{C_1}{A_1} \\&1&-\frac{B_2}{A_2}\\&&1\end{matrix}\right). \]
Note that we put our decompositions in the form $bcwb'$ with $b' \in \bar{U}_w(\Bbb{Q})$, as opposed to \cite{BFG}, who put their decompositions in the form $bwcb'$ with $b \in \bar{U}_{w^{-1}}(\Bbb{Q})$.
The decompositons are equivalent, but the former is more standardized.

Restricting to a fundamental domain for the action of $\bar{U}_{w_6}(\Bbb{Z}) = U(\Bbb{Z})$, we may write the sum \eqref{kloabs} for $w=w_6$ as 
\[ S_{w_6}(\psi_{n_1,n_2},\psi_{m_1,m_2};(A_1,A_2)) = \sum_{B_1,C_1,B_2,C_2} e\left(n_1 \frac{Y_1 A_2-Z_1 B_2}{A_1}+n_2 \frac{Z_2 B_1 - Y_2 A_1}{A_2}-m_1 \frac{B_2}{A_2}+m_2 \frac{B_1}{A_1}\right), \]
where the sum is taken over $0 \le B_1, C_1 < A_1$ and $0 \le B_2 < A_2$ subject to \eqref{eq:PlukerConds} and \eqref{eq:AuxPluker}.
Now \cite[Lemma 4.1]{BFG} shows that the summand is independent of the choice of auxiliary parameters.
Note that the compatibility condition is trivially true for the long element as $U_{w_6}(\Bbb{R}) = \{I\}$.

We conclude the long element analysis by mentioning that the proofs of \cite[Lemmas 4.1, 4.2]{BFG} show that the sum is well-defined if we replace the summation conditions with their modular equivalents
\begin{align}
\label{eq:LongElementConds}
\begin{array}{l}
	B_1, C_1 \pmod{A_1}, \qquad B_2 \pmod{A_2}, \qquad N\mid A_1, \qquad N \mid  B_1\\
	A_1 C_2+B_1B_2+C_1 A_2\equiv0 \pmod{A_1 A_2}, \\
	(A_1,B_1,C_1)=1, \qquad (A_2,B_2,C_2)=1, \\
	Z_2 C_2+Y_2 B_2=1\pmod{A_2}, \qquad Z_1 C_1+Y_1 B_1=1 \pmod{A_1},
\end{array}
\end{align}
and the sum is empty unless $N\mid A_2$.
This matches the previous definition with
\begin{equation}\label{klo6conv}
S_{w_6}(\psi_{n_1,n_2},\psi_{m_1,m_2};(A_1,A_2)) = \delta_{N\mid A_1} S^{(N)}(m_2,-m_1,n_1,-n_2;A_1, A_2).
\end{equation}
As noted in the proof of \cite[Theorem 5.1]{BFG}, replacing $m_1,n_2,B_2,C_2,C_1,Y_2,Z_2,Z_1$ by their negatives leaves the sum invariant, so we may drop the negatives on $m_1$ and $n_2$.

Elements of the $w_5$ cell necessarily have $A_1 = 0$ and $B_1,A_2 \ne 0$, so that $B_1 \mid  A_2$.
With the Pl\"ucker and auxiliary coordinates as before, for $\gamma$ having $A_1=0$ we have
\[ \gamma = \left(\begin{matrix}1&\frac{Z_2 B_1}{A_2}&\frac{Y_2}{B_1} \\&1&\frac{Z_1 C_2-X_1 A_2}{B_1}\\&&1\end{matrix}\right) \left(\begin{matrix}\frac{1}{A_2}\\&\frac{A_2}{B_1}\\&&B_1\end{matrix}\right) w_5 \left(\begin{matrix}1&0&-\frac{C_2}{A_2} \\&1&\frac{C_1}{B_1}\\&&1\end{matrix}\right). \]
The compatibility condition becomes $n_1 A_2=m_2 B_1^2$.
Now the conditions \eqref{eq:PlukerConds} and \eqref{eq:AuxPluker} simplify to
\begin{align*}
\begin{array}{l}
	(B_1,C_1)=1, \qquad (A_2/B_1,C_2)=1, \qquad B_2 = -C_1 \frac{A_2}{B_1}, \qquad N\mid B_1, \\
	Z_1 C_1 \equiv 1 \pmod{B_1}, \qquad Z_2 C_2 \equiv 1 \pmod{A_2/B_1},
\end{array}
\end{align*}
and the space $\bar{U}_{w_5}(\Bbb{Z}) \subset U(\Bbb{Z})$ is defined by $x_2=0$, so we may write the Kloosterman sum as
\[ S_{w_5}(\psi_{n_1,n_2},\psi_{m_1,m_2};(B_1,A_2)) = \delta_{\substack{n_1 A_2=m_2 B_1^2 \\ N\mid B_1\mid A_2}} \tilde{S}(m_1,n_1,n_2;B_1,A_2). \] 
Unlike the long element case, no extra work is needed to justify our use of the modular summation conditions.

For the $w_4$ cell, we have $A_2=0$ and $A_1,B_2\ne 0$, so that $B_2\mid A_1$, and the compatibility condition is $n_2 A_1=m_1 B_2^2$.
With the Pl\"ucker and auxiliary coordinates as before, we have
\[ \gamma = \left(\begin{matrix}1&\frac{X_2 A_1-Z_2 C_1}{B_2}&\frac{Z_2}{A_1} \\&1&-\frac{Z_1 B_2}{A_1}\\&&1\end{matrix}\right) \left(\begin{matrix}\frac{1}{B_2}\\&\frac{B_2}{A_1}\\&&A_1\end{matrix}\right) w_4 \left(\begin{matrix}1&-\frac{C_2}{B_2}&\frac{C_1}{A_1} \\&1&0\\&&1\end{matrix}\right). \]
Now the conditions \eqref{eq:PlukerConds} and \eqref{eq:AuxPluker} simplify to
\begin{align*}
\begin{array}{l}
	(A_1/B_2,C_1)=1, \qquad (B_2,C_2)=1, \qquad B_1 = -C_2 \frac{A_1}{B_2}, \qquad N\mid A_1, \qquad N\mid B_1, \\
	Z_1 C_1 \equiv 1 \pmod{A_1/B_2}, \qquad Z_2 C_2 \equiv 1 \pmod{B_2},
\end{array}
\end{align*}
but we may take this one step further: 
we have $B_2\frac{B_1}{N} = -C_2 \frac{A_1}{N}$ and the condition $(B_2,C_2)=1$ implies $B_2 \mid  \frac{A_1}{N}$.
Conversely, the condition $N B_2 \mid  A_1$ implies $N \mid  (-C_1 \frac{A_1}{B_2})=B_1$, so we may write the Kloosterman sum as
\[ S_{w_4}(\psi_{n_1,n_2},\psi_{m_1,m_2};(A_1,B_2)) = \delta_{\substack{n_2 A_1=m_1 B_2^2 \\ N B_2\mid A_1}} \tilde{S}(-m_2,-n_2,-n_1;B_2,A_1). \] 
Note that changing the sign of both $n_1$ and $n_2$ leaves $\tilde{S}$ invariant by $C_2 \mapsto -C_2$.

\section{Weight functions in the Kuznetsov formula}

In order to  use the Kuznetsov formula for a spectral average, we need a function $F$ such that $|\langle F, \tilde{W}_{\mu}\rangle|^2$ appearing on the left hand side of \eqref{kuzform} is bounded away from zero for $\mu \in \Omega \subseteq \mathfrak{a}_{\Bbb{C}}^{\ast}$. For our purposes the following slightly weaker statement suffices.

\begin{lemma}\label{omega} For a fixed compact $\Omega \subseteq \mathfrak{a}_{\Bbb{C}}^{\ast}$ there is a finite collection of smooth compactly supported functions $F_1, \ldots, F_J$ such that $\sum_j |\langle F_j, \tilde{W}_{\mu}\rangle|^2 \gg 1$ for $\mu \in \Omega.$ 
\end{lemma}

\textbf{Proof.} This follows from a simple compactness argument: for each $\mu \in \Omega$ choose an open set $S_{\mu} \subseteq \Bbb{R}_{>0}^2$ such that $\Re \tilde{W}_{\mu}(y) \not= 0$ for all $y \in S_{\mu}$ or  $\Im \tilde{W}_{\mu}(y)\not= 0$  for all $y \in S_{\mu}$. This is possible by continuity of $\tilde{W}_{\mu}(y)$ in $y$. By continuity in $\mu$, we can choose open neighbourhoods $U_{\mu}$ about $\mu$ such that  $\Re \tilde{W}_{\mu^{\ast}}(y) \not= 0$ for all $y \in S_{\mu}$ and all $\mu^{\ast} \in U_{\mu}$ or  $\Im \tilde{W}_{\mu}(y)\not= 0$  for all $y \in S_{\mu}$ and all $\mu^{\ast} \in U_{\mu}$. By compactness we pick a finite collection  of such neighbourhoods $U_{\mu_1}, \ldots, U_{\mu_J}$ covering $\Omega$, and define the corresponding $F_j$ to be real-valued functions with support on $S_{\mu_j}$ and non-vanishing on the interior $\mathring{S}_{\mu_j}$.    \\

For the proof of Theorem \ref{thm4} we will  need a function that blows up on the exceptional spectrum. We recall that by unitarity the exceptional spectrum is parametrized by 
\begin{equation}\label{except}
\mu = (\rho + i\gamma, - \rho + i \gamma, -2i\gamma)
\end{equation}
for $\gamma \in \Bbb{R}$, $\rho \in [-1/2, 1/2]$ (by the Jacquet-Shalika bounds) 
   and its translates under the Weyl group. For a fixed smooth compactly supported function $F$ and two parameters $X_1, X_2 \geq 1$ let
 \begin{equation}\label{shift}  
   F^{(X_1, X_2)}(y_1, y_2) := F(X_1y_1, X_2 y_2)
   \end{equation} 
   so that $F = F^{(1, 1)}$. 

 \begin{lemma}\label{X} Fix $\Omega  \subseteq \mathfrak{a}_{\Bbb{C}}^{\ast}$, and let $X_1, X_2 \geq 1$, $\varepsilon > 0$.  Assume that $F$ is non-negative and supported  in a (depending on $\Omega$) sufficiently small neighbourhood about $(1, 1)$ and that $X_1, X_2$ are sufficiently large.  Then for exceptional $\mu \in \Omega$ of the form \eqref{except} with $|\rho| \geq \varepsilon$ we have
$\langle F^{(X_1, X_2)}, \tilde{W}_{\mu} \rangle \asymp (X_1X_2)^{1+|\rho|}$.   
\end{lemma}

 \textbf{Proof.}  We have by \eqref{whit} that
\begin{displaymath}
\begin{split}  
 &   \langle F^{(X_1, X_2)}, \tilde{W}_{\mu} \rangle\\
    & =   \int_{(1)} \int_{(1)} \frac{\cosh(\frac{3}{2} \pi \gamma)\prod_{j=1}^3 \Gamma(\frac{1}{2}(s_1 + \mu_j)) \prod_{j=1}^3 \Gamma(\frac{1}{2}(s_2 - \mu_j)) }{4\pi^{s_1+s_2 } \Gamma(\frac{1}{2}(s_1+s_2))} \mathcal{F}(-1-s_1, -1-s_2) X_1^{1+s_1} X_2^{1+s_2} \frac{ds_1\, ds_2 }{(2\pi i)^2}
    \end{split}
    \end{displaymath}
where $\mathcal{F}$ is the double Mellin transform of $F$, an entire function in both variables.  If  without loss of generality $\rho > 0$ (note that in particular $\mu_1, \mu_2, \mu_3$ are sufficiently distinct), we shift contours to the left and obtain
$$  \langle F^{(X_1, X_2)}, \tilde{W}_{\mu} \rangle = c_{\mu} \mathcal{F}(-1-\rho + i\gamma, - 1 - \rho - i\gamma) (X_1^{1+\rho - i\gamma}  + O(X_1)) (X_2^{1+\rho + i\gamma} + O(X_2))$$
for some constant $c_{\mu} \not= 0$. If $F$ is non-negative and supported in a sufficiently small   neighbourhood about $(1, 1)$, then $\mathcal{F}(-1-\rho + i\gamma, - 1 - \rho - i\gamma) \not= 0$ for all $\mu \in \Omega$. This proves the lemma.\\

Next we provide bounds for the functions  $\tilde{\mathcal{J}}_{\epsilon;F}(A)$ and $ \mathcal{J}_{\epsilon;F} (A_1, A_2)$ defined in \eqref{test1} and \eqref{test2}. Here $F$ will always be a fixed compactly supported function and all implied constants may depend on $F$. For bounds in the case of certain highly oscillating functions $F$ see \cite[Proposition 5]{Bl}. We define $F^{(X_1, X_2)}$ as in \eqref{shift}. The following basic bound suffices in many cases. 

\begin{lemma}\label{J1} Let $X_1, X_2 \geq 1$. \\
 {\rm (a)} We have  $\tilde{\mathcal{J}}_{\epsilon; F^{(X_1, X_2)}}(A) = 0$  unless $A \gg X_1^{-3/2} + X_2^{-3/2}$, in which case $$  \tilde{\mathcal{J}}_{\epsilon; F^{(X_1, X_2)}}(A) \ll (X_1X_2)^2. $$

\noindent {\rm (b)} We have  $\mathcal{J}_{\epsilon;F^{(X_1, X_2)}}(A_1, A_2) = 0$  unless $\min(A_1 A_2^2, A_2A_1^2)  \gg (X_1X_2)^{-3/2}$, in which case $$\frac{d^i}{dA_1^i}  \frac{d^j}{  dA_2^j} \mathcal{J}_{\epsilon; F^{(X_1, X_2)}}(A_1, A_2) \ll_{i, j}  (X_1X_2)^2(A_1A_2)^{\varepsilon} \left(A_2^{2/3} A_1^{1/3}\right)^i\left(A_1^{2/3} A_2^{1/3}\right)^j$$
for all $i, j \in \Bbb{N}_0 = \{0,1, 2, \ldots\}$. 
\end{lemma}

\textbf{Remark:} Except for one instance in the proof of Theorem \ref{thm4} we will always apply this lemma with $X_1 = X_2 = 1$, so for  most of the paper the variables $X_1, X_2$ can be ignored. \\

\textbf{Proof.} (a) This is straightforward from the definition and uses only trivial bounds, noting that the support of $F^{(X_1, X_2)}$ restricts the variables to
$$y_1 \asymp (X_1A)^{-1}, \quad y_2 \asymp X_2^{-1}, \quad 1 + x_1^2 \asymp A^{4/3}X_1^2, \quad 1 + x_1^2 + x_2^2 \asymp A^{8/3}(X_1X_2)^2.$$
This forces $A \gg X_1^{-3/2} + X_2^{-3/2}$. The upper bound follows now from trivial estimates\footnote{This corrects an error \cite[(8.7)]{Bl} where $X_1^2 X_2$ should be replaced with $X_1^2X_2^2$.}.

(b)  The support of $F^{(X_1, X_2)}$ restricts the variables to
\begin{displaymath}
\begin{split}
& y_1 \asymp (X_1A_1)^{-1}, \quad y_2 \asymp (X_2A_2)^{-1}, \\
& (x_1x_2 - x_3)^2 + x_1^2 + 1 =: \xi_1  \asymp \Xi_1 := A_2^{4/3} A_1^{8/3}(X_1X_2)^2,
\quad  x_3^2 + x_2^2 + 1 =:\xi_2 \asymp \Xi_2 := A_1^{4/3} A_2^{8/3}(X_1X_2)^2
\end{split}
\end{displaymath}
which implies that both $A_1A_2^2$ and $A_1^2 A_2$ must be at least of order $(X_1X_2)^{-3/2}$. 
We recall from \cite[Lemma 4]{Bl} that
\begin{equation}\label{lemma4}
\int_{\substack{\xi_1 \asymp \Xi_1\\ \xi_2 \asymp \Xi_2}} dx_1\, dx_2\, dx_3 \ll (\Xi_1\Xi_2)^{1/2 + \varepsilon} =(A_1A_2X_1X_2)^{2+\varepsilon}
\end{equation}
for $\Xi_1, \Xi_2 \gg 1$. For the derivatives we differentiate under the integral sign and estimate  trivially, see also \cite[(8.16)]{Bl}. \\

 
 For one application we need a more refined estimate of a certain $6$-fold Fourier transform involving $\mathcal{J}_{\epsilon;F}$. 
 
 \begin{lemma}\label{stat2}  Let $W : (0, \infty)^6 \rightarrow\Bbb{C}$ be a fixed smooth compactly supported function. Let $A_1, A_2 > 0$ and define $A := \exp(\max (|\log A_1|, |\log A_2|))$. Let $P \geq 1$, and let $\alpha_1, \alpha_2, \beta_1, \beta_2, \gamma_1, \gamma_2 \in \Bbb{R}$ be   such that $\min(|\alpha_1|, |\alpha_2|, |\beta_1|, |\beta_2|, |\gamma_1|, |\gamma_2| ) \leq P$.   Then the six-fold Fourier transform
\begin{displaymath}
\begin{split}
\widehat{\mathcal{J}} :=   \int_{\Bbb{R}^6} \mathcal{J}_{\epsilon; F}(&A_1 \sqrt{t_1u_1v_1}, A_2\sqrt{t_2u_2v_2}) W(t_1, t_2, u_1, u_2, v_1, v_2)\\
 &\times e(-t_1\alpha_1 - t_2\alpha_2 - u_1\beta_1-u_2\beta_2 - v_1\gamma_1-v_2\gamma_2) dt_1\, dt_2\, du_1\, du_2\, dv_1\, dv_2
 \end{split}
 \end{displaymath}
is bounded by 
\begin{equation}\label{new1}
O_C\left((PA)^{\varepsilon}( P^2 \max\left(A_2^{-2/3} A_1^{-4/3}, A_1^{-2/3} A_2^{-4/3}\right)+ P^{-C})\right)
\end{equation}
for any constant $C > 0$. In addition, it is bounded by
\begin{equation}\label{new2}
A^{\varepsilon}\max(|\alpha_1|,   |\beta_1|,   |\gamma_1|  )^{-1/2} \max(|\alpha_2|,   |\beta_2|,   |\gamma_2|  )^{-1/2},
\end{equation}
as long as both maxima are non-zero.
\end{lemma}

\textbf{Proof.}  We recall the notation
$$ (x_1x_2 - x_3)^2 + x_1^2 + 1 =: \xi_1  \asymp \Xi_1 := A_2^{4/3} A_1^{8/3}, 
\quad  x_3^2 + x_2^2 + 1 =:\xi_2 \asymp \Xi_2 := A_1^{4/3} A_2^{8/3}
$$
from the previous proof. We will frequently use \eqref{lemma4} with $X_1 = X_2 = 1$. 
We also write $x_4 = x_1x_2 - x_3$ and 
$$\eta_1 = x_2x_4 + x_1, \quad \eta_2 = x_1x_3 + x_2.$$
We   express  $ \mathcal{J}_{\epsilon; F}$ by its defining 5-fold integral \eqref{test2} and write
$$\mathcal{J}_{\epsilon; F}(A_1, A_2) =  \frac{1}{(A_1A_2)^2}\int_{\substack{\xi_1 \asymp \Xi_1\\ \xi_2 \asymp \Xi_2}} \mathcal{K}(A_1, A_2; x_1, x_2, x_3) dx_1\, dx_2\, dx_3$$
 where $\mathcal{K}$  is the double $y_1, y_2$-integral, i.e.
 \begin{equation}\label{added-int}
 \begin{split}
 \mathcal{K}(A_1, A_2; x_1, x_2, x_3)  = \int_0^{\infty}\int_0^{\infty}  &e\left(-\epsilon_1 A_1x_1y_1 - \epsilon_2 A_2x_2y_2 -\frac{A_2\eta_2}{y_2\xi_2}   -\frac{A_1\eta_1}{y_1\xi_1}   \right)\\
 & \times \overline{F(A_1y_1, A_2 y_2)}   F\left(\frac{A_2\xi_1^{1/2}}{y_2\xi_2} , \frac{A_1\xi_2^{1/2}}{y_1\xi_1}  \right)   \frac{dy_1\, dy_2}{y_1y_2}. 
\end{split}
 \end{equation}
 
 We start with the proof of \eqref{new1}.  Suppose that $|\alpha_1|$ is the smallest of the variables (possibly $|\alpha_1| = 0$). Choose a sufficiently large constant $c_2$ and a sufficiently large constant $c_1 > c_2$. We split the $x_1, x_2, x_3$-integration in four pieces 
$$ {\rm (i)}\,\,  |x_1|, A_1^2 |\eta_1|/\xi_1 \leq c_1 P, \quad {\rm (ii)}\,\, |x_1| \leq c_2P, A_1^2 |\eta_1|/\xi_1 \geq c_1 P, \quad {\rm (iii)} \,\,|x_1| \geq c_1P, A_1^2 |\eta_1|/\xi_1 \leq c_2P$$
and the remaining portion (iv), which is contained in $
 |x_1|, A_1^2 |\eta_1|/\xi_1 \geq c_2 P$. The conditions (i) imply $|x_1| \ll P$, $x_2x_4 \ll P A_2^{4/3} A_1^{2/3}$ (note that we may assume by Lemma \ref{J1} that this is  $\gg P$). The area of this region is
 $$\ll P \int_{\substack{x_2x_4 \ll P A_2^{4/3} A_1^{2/3}\\ x_2, x_4 \ll A^{O(1)}}} dx_2\, dx_4 \ll   P^2A_2^{4/3} A_1^{2/3} (AP)^{\varepsilon}.$$
As  $\mathcal{K}(A_1, A_2; x_1, x_2, x_3) \ll 1$, the total contribution of this case to $\widehat{\mathcal{J}}$ is $P^2A_2^{-2/3} A_1^{-4/3} (AP)^{\varepsilon}$. 

To deal with the region (ii), we note that the phase in the $y_1$-integral in \eqref{added-int} is given by 
$$e\left( - \epsilon_1A_1x_1y_1 - \frac{A_1\eta_1}{y_1\xi_1}\right).$$
 If $c_1$ is sufficiently large compared to $c_2$ (or if $\epsilon_1 x_1$ and $\eta$ have different signs), the phase has no stationary point, and after sufficiently many integrations by parts,  using for instance \cite[Lemma 8.1]{BKY} with
$$X = A_1, \quad U = Q = \frac{1}{A_1}, \quad Y = P, \quad R = A_1Y,$$
we bound after trivial estimation in all other variables this portion of $\widehat{\mathcal{J}}$ by $ \ll A^{\varepsilon}P^{-C}.$ The same argument works for the region (iii).

 In order to analyze the region (iv), we consider the expression \eqref{added-int} in more detail, first without any restrictions on the $x$-variables. We could run a careful stationary phase argument as in \cite[Proposition 8.2]{BKY}, but we can also proceed in a completely elementary way. Applying the stationary phase method only on a formal basis shows that the oscillation of the $y_1$-integral is given by $$e\Bigl(- \text{sgn}(\eta_1) \frac{2 \sqrt{|x_1\eta_1|} A_1}{\sqrt{\xi_1}}\Bigr), $$
coming from the stationary point at $y_1 = (\eta_1/(\epsilon_1 x_1\xi_1))^{1/2}$. 
With this in mind let us define  
\begin{displaymath}
\begin{split}
\tilde{\mathcal{K}}(A_1, A_2; x_1, x_2, x_3) & := e\left( \text{sgn}(\eta_1)   \frac{2 \sqrt{|x_1\eta_1|} A_1}{\sqrt{\xi_1}}+ \text{sgn}(\eta_2)  \frac{2 \sqrt{|x_2\eta_2|} A_2}{\sqrt{\xi_2}}\right)\mathcal{K}(A_1, A_2; x_1, x_2, x_3)\\
& = 
  \int_0^{\infty}\int_0^{\infty}   e\Bigl(g(A_1, y_1) + h(A_2, y_2)\Bigr)
  \overline{F(A_1y_1, A_2 y_2)}   F\left(\frac{A_2\xi_1^{1/2}}{y_2\xi_2} , \frac{A_1\xi_2^{1/2}}{y_1\xi_1}  \right)   \frac{dy_1\, dy_2}{y_1y_2}
\end{split}
\end{displaymath}
with
$$g(A_1, y_1) = g_{\epsilon_1 x_1, \eta_1, \xi_1}(A_1, y_1) =  - \epsilon_1 A_1x_1y_1 - \frac{A_1\eta_1}{y_1\xi_1} +  \text{sgn}(\eta_1) \frac{ 2\sqrt{|x_1\eta_1|}A_1  }{\sqrt{\xi_1}}$$
and
$$h(A_2, y_2) = h_{\epsilon_2 x_2, \eta_2, \xi_2}(A_2, y_2) =  - \epsilon_2 A_2x_2y_2 - \frac{A_2\eta_2}{y_2\xi_2} +  \text{sgn}(\eta_2) \frac{ 2\sqrt{|x_2\eta_2|}A_2 }{\sqrt{\xi_2}}.$$
We show the uniform bound 
\begin{equation}\label{eq-in-lemma4}\frac{\partial^i}{\partial A_1^i}  \frac{\partial^j}{\partial A_2^j}\tilde{\mathcal{K}}(A_1, A_2; x_1, x_2, x_3) \ll_{i, j} A_1^{-i}A_2^{-j}. 
\end{equation}
Indeed, one checks by direct computation that
$$  \frac{\frac{\partial}{\partial A_1} g(A_1, y_1)}{\frac{\partial}{\partial y_1} g(A_1, y_1)}   = \pm 
\frac{y_1}{A_1} \cdot  \frac{\sqrt{|\eta_1|} - \sqrt{|x_1|\xi_1}y_1}{\sqrt{|\eta_1|} + \sqrt{|x_1|\xi_1}y_1} $$
so that
$$
 \frac{\partial^i}{\partial y^i_1}\frac{\partial^j}{\partial A^j_1} \Bigl( \frac{\frac{\partial}{\partial A_1} g(A_1, y_1)}{\frac{\partial}{\partial y_1} g(A_1, y_1)}\Bigr)   \ll_{i, j}  \frac{y_1}{A_1} y_1^{-i}A_1^{-j}$$
 for $i, j \in \Bbb{N}_0$. 
Hence combining each differentiation with respect to $A_1$ with an integration by parts in $y_1$, we obtain the desired bound \eqref{eq-in-lemma4} in $A_1$, and the bound in $A_2$ follows similarly. 

Having proved \eqref{eq-in-lemma4}, we return to the estimation of $\widehat{\mathcal{J}}$. The phase of the Fourier integral in question is given by 
\begin{equation}\label{phase-new2}
e\left(\pm \frac{ 2\sqrt{|x_1\eta_1|} A_1 \sqrt{t_1u_1v_1}}{\sqrt{\xi_1}}\pm  \frac{ 2\sqrt{|x_2\eta_2|} A_2 \sqrt{t_2u_2v_2}}{\sqrt{\xi_2}} -t_1\alpha_1 - t_2\alpha_2 - u_1\beta_1-u_2\beta_2 - v_1\gamma_1-v_2\gamma_2\right),
\end{equation}
which needs to be integrated against the non-oscillating functions
$$ \tilde{\mathcal{K}}(A_1\sqrt{t_1u_1v_1}, A_2\sqrt{t_2u_2v_2}; x_1, x_2, x_3) W(t_1, t_2, u_1, u_2, v_1, v_2)$$
with respect to $x_1, x_2, x_3$  and $t_1, t_2, u_1, u_2, v_1, v_2$. If we are in the region (iv), then in particular $
 |x_1|, A_1^2 |\eta_1|/\xi_1 \geq c_2 P$, as mentioned above. Since $|\alpha_1| \leq P$,   the phase has no stationary point if $c_2$ is sufficiently large, and by repeated partial integration in any of the variables $t_1$ we obtain again the bound $A^{\varepsilon} P^{-C}$. This completes the proof of \eqref{new1} if $|\alpha_1|$ is minimal. If any of the other variables is minimal, we can run the same argument, possibly with interchanged indices. 
 
For the bound \eqref{new2} we return to \eqref{phase-new2} for an arbitrary choice of  $x_1, x_2, x_3$. The simple  stationary phase type bound
 $$\int e(at+b \sqrt{t}) W(t) dt \ll |a|^{-1/2}, \quad a \not= 0$$
for a fixed smooth function $W$ with compact support in $(0, \infty)$ applied twice, followed by trivial estimations, yields readily the bound \eqref{new2}.  This completes the proof.

\section{Kloosterman sums}

In this section we collect some  results about the  Kloosterman sums defined in \eqref{klo45} and  \eqref{klo6}.   We start with useful upper bounds. 
 
\begin{lemma}\label{lem1} Let $N, D_1, D_2 \in \Bbb{N}$, $m_1, m_2, n_1, n_2 \in \Bbb{Z} \setminus \{0\}$.  We have
$$\tilde{S}(m_1, m_2, n_1; D_1, D_2) \ll \Bigl((n_2, D_2/D_1)D_1^2, (m_1, n_1, D_1)D_2\Bigr) (D_1D_2)^{\varepsilon}$$
 and
 \begin{equation}\label{Stevens}
S^{(N)}(m_1, m_2, n_1, n_2;  D_1, D_2) \ll (D_1D_2)^{1/2 + \varepsilon}\bigl((D_1, D_2)(m_1n_1, [D_1, D_2])(m_2n_2, [D_1, D_2])\bigr)^{1/2}
\end{equation}
for any $\varepsilon > 0$. 
\end{lemma} 
 
\textbf{Proof.} The bound for $\tilde{S}$ is Larsen's bound \cite[Appendix]{BFG}. The bound for $S^{(N)}$ is Stevens' bound \cite[Theorem 5.1]{St} in its uniform version given in \cite[p.\ 39]{Bu}. Note that for the level $N$ Kloosterman sum, only those $S_{a, b}(n, \psi, \psi')$ (in the notation of \cite[Section 5]{St}) contribute to the Kloosterman sum where $a \leq s-k$, $b \leq r$ with $s, r \geq k$ whenever $p^k \parallel N$ (cf.\ also \cite[Remark 2.5]{DF}). In particular, Stevens' bound holds a fortiori for level $N$ Kloosterman sums.\\ 

As in \cite[Lemma 3]{Bl}\footnote{Notice that  in \cite[Lemma 3]{Bl} the indices should be exchanged and read as in \eqref{Stevens} above as a consequence of Remark 1 after Theorem \ref{kuz}.} we conclude from \eqref{Stevens} that 
\begin{equation}\label{average}
\begin{split}
&\sum_{\substack{N \mid D_1 \leq X_1 \\ N \mid D_2 \leq X_2}}  |S^{(N)}(m_1, m_2, n_1, n_2; D_1, D_2) | \\
&\ll (X_1X_2)^{1/2+\varepsilon}N^{1/2} \sum_{\substack{\delta \delta_1 \leq X_1/N\\ \delta \delta_2 \leq X_2/N}}\delta^{1/2} \bigl((m_1n_1, \delta_1)(m_1n_1, \delta_2) (m_1n_2, \delta) (m_2n_2, \delta_1) (m_2n_2, \delta_2) (m_2n_2, \delta)\bigr)^{1/2}   \\
&\ll \frac{(X_1X_2)^{3/2+\varepsilon}(m_1n_2 m_2n_1)^{\varepsilon}}{N^{3/2}} \sum_{\delta \leq X_1} \frac{(m_1n_1, \delta)^{1/2}(m_2n_2, \delta)^{1/2}}{\delta^{3/2}} \ll \frac{(X_1X_2)^{3/2+\varepsilon}(m_1n_2 m_2n_1)^{\varepsilon}}{N^{3/2}} 
\end{split}
\end{equation}
if $(m_1m_2n_1n_2, N) = 1$.

\begin{lemma}\label{prop1} For $N \in \Bbb{N}$ the following  holds.\\
{\rm (a)} The sum $S^{(N)}(m_1, m_2, n_1, n_2; D_1, D_2)$ depends only on $m_j, n_j$ modulo $D_j$ for $j=1, 2$. \\
{\rm (b)} If   $(t_1t_2, u_1u_2) = 1$ for $j=1, 2$ and $N \mid t_1u_1, t_2u_2$, then
\begin{displaymath}
\begin{split}
&S^{ (N)}(m_1, m_2, n_1, n_2; t_1u_1, t_2u_2) \\
&= S^{(\text{{\rm gcd}}(N, t_1))}(\bar{u}_1^2 u_2 m_1, \bar{u}_2^2 u_1 m_2, n_1, n_2; t_1,t_2) S^{({\rm gcd}(N, u_1))}(\bar{t}_1^2t_2m_1, \bar{t}_2^2t_1m_2, n_1, n_2; u_1, u_2).
\end{split}
\end{displaymath}
{\rm (c)} Let $N$ be prime and let $r_q(n)$ denote the Ramanujan sum. Then
$$S^{(N)}(m_1, m_2, n_1, n_2; N, N)= N-1 + r_N(n_1) r_N(m_2) = \begin{cases}N(N-1), & N \mid n_1, N \mid m_2,\\ N , & N \nmid n_1m_2,\\ 0 \, & \text{else.} \end{cases}$$ 
\end{lemma}

\textbf{Proof.} This is proved as in   \cite[Section 4]{BFG}, cf.\ Properties 4.6,  4.7 and 4.10, respectively. \\
 Part (a) is trivial. \\
For part (b) we observe that the   assumptions $(t_1, u_1) = 1$ and $N \mid t_1u_1$ imply that 
\begin{equation}\label{gcd-new}
\big((N, t_1)u_1, (N, u_1)t_1\big) = N.
\end{equation}
We now follow verbatim the proof of \cite[Property 4.7]{BFG}. Given two sets of summation variables $B_j, C_j$ and $B_j', C_j'$ such that $(B_j, C_j, t_j) = (B_j', C_j', u_j) = 1$ for $j=1, 2$, $(N, t_1) \mid B_1$, $(N, u_1) \mid B_1'$, $t_1t_2 \mid t_1C_2 + B_1B_2 + C_1 t_2$ and $u_1u_2 \mid u_1C'_2 + B'_1B'_2 + C'_1 u_2$ we choose $r, r' \in \Bbb{Z}$ with $rt_1t_2 + r' u_1u_2 = 1$. We   define new variables 
\begin{displaymath}
\begin{split}
&  d_1 = t_1u_1, \quad b_1 = r' u_1u_2B_1 + r t_1t_2 B_1', \quad c_1 = (r')^2 u_1^2 u_2 C_1 + r^2 t_1^2 t_2 C_1',\\
 & d_2 = t_2u_2, \quad b_2 = r' u_1u_2B_2 + r t_1t_2 B_2', \quad c_2 = (r')^2 u_1 u_2^2 C_2 + r^2 t_1 t_2^2 C_2',\\
\end{split}
\end{displaymath}
and observe that $b_1$ runs through all numbers modulo $d_1 = t_1u_1$ that are divisible by \eqref{gcd-new}, which is the desired extra divisibility condition $N \mid b_1$ for  $S^{ (N)}(m_1, m_2, n_1, n_2; t_1u_1, t_2u_2) $.  Now we continue verbatim as in \cite{BFG}.\\
 To prove (c) we observe that   $S^{(N)}(m_1, m_2, n_1, n_2; N, N) $  for $N$ prime equals
 \begin{displaymath}
\begin{split}
 &   \sum_{\substack{C_1, B_2, C_2 \, ({\rm mod }\,  N)\\  C_2   +  C_1 \equiv 0 \, ({\rm mod }\,  N)\\ (  C_1, N) = (B_2, C_2, N)   = 1 }} e\left(\frac{  -n_1 \bar{C}_1 B_2}{N} + \frac{m_2B_2  }{N}\right)=   \sum_{\substack{B, C \, ({\rm mod }\, N)\\ (C, N) = 1}} e\left(\frac{  -n_1 C  B}{N} + \frac{m_2B }{N}\right)\\
 & =  N-1 + \sum_{\substack{B, C\, (\text{mod } N)\\ (BC, N) = 1}} e\left(\frac{  -n_1 C  B}{N} + \frac{m_2B }{N}\right)  = N-1 + r_N(n_1) r_N(m_2).    \end{split}
\end{displaymath}
This completes the proof of the lemma. \\

\textbf{Remark:} For completeness we also state the following two properties of the level $N$ Kloosterman sums that can be proved as in Property 4.3 and Property  4.4+4.5 of \cite{BFG}:\\
(d) For $(D_1D_2, a b) = 1$ we have
$S^{( N)}(am_1, bm_2, n_1, n_2; D_1, D_2)  = S^{(N)}(m_1, m_2, an_1, bn_2; D_1, D_2).$\\
(e) We have $S^{( N)}(m_1, m_2, n_1, n_2; D_1, D_2) = S^{( N)}(n_2, n_1, m_2, m_1; D_2, D_1).$\\
Note, however, that an analogue of Property 4.4 or Property 4.5 alone does not exist due to asymmetry of the summation condition $N \mid B_1$. We do not need the statements (d) and (e) in this paper. \\

For later purposes we study a certain 6-fold Fourier transform of the long Weyl element Kloosterman sum. Let $d, D_1, D_2 \in \Bbb{N}$,   and let $a \in (\Bbb{Z}/D_1Z)^{\ast}$ and $b \in (\Bbb{Z}/D_2\Bbb{Z})^{\ast}$. 
For integers $x_1, x_2, y_1, y_2, z_1, z_2$ we define 
\begin{equation}\label{untwist}
\begin{split}
\widehat{S}_{a, b, d}(x_1, x_2, y_1, y_2, z_1, z_2; D_1, D_2) := & \frac{1}{D_1^3D_2^3}\sum_{\substack{n_1, m_1, l_1 \, (\text{mod } D_1) \\ n_2, m_2, l_2 \, (\text{mod } D_2)}}   S^{(1)}(am_1d, bn_2l_2, n_1l_1, m_2 d;D_1, D_2) \\
& \times e\left(-\frac{n_1x_1 + m_1y_1+l_1z_1}{D_1}\right) e\left(-\frac{n_2x_2 + m_2y_2+l_2z_2}{D_2}\right).
 \end{split}
 \end{equation}
 This is the non-archimedean analogue of the function studied in Lemma \ref{stat2}. 
 
We also need a twisted version. 
Let $\chi$ be a primitive character modulo a prime $p$ such that $(d, p) = 1$. Assume that $p^3 \mid D_1$ and $p^3 \mid D_2$. Then we define
\begin{equation}\label{shat}
\begin{split}
\widehat{S}^{\chi}_{a, b, d}&(x_1, x_2, y_1, y_2, z_1, z_2; D_1, D_2) := \frac{1}{D_1^3D_2^3}\sum_{\substack{ n_1, m_1, l_1 \, (\text{mod } D_1)\\n_2, m_2, l_2 \, (\text{mod } D_2)}} \bar{\chi}(n_1l_1m_2)\chi(m_1n_2l_2) \\
&\times S^{(p^3)}(am_1d, bn_2l_2, n_1l_1, m_2 d;D_1, D_2) e\left(-\frac{n_1x_1 + m_1y_1+l_1z_1}{D_1}\right) e\left(-\frac{n_2x_2 + m_2y_2+l_2z_2}{D_2}\right).
 \end{split}
 \end{equation}
 By  the Chinese remainder theorem and Lemma \ref{prop1}(b) we have the following multiplicativity formulae
 \begin{equation}\label{mult1}
 \begin{split}
 \widehat{S}^{\chi}_{a, b, d}(x_1, x_2, y_1, y_2, z_1, z_2;  t_1p^{\alpha_1}, t_2p^{\alpha_2}) & = \widehat{S}^{\chi}_{a, b, d}(\bar{t}_1x_1, t_2\bar{t}_1x_2, t_1\bar{t}_2y_1,  \bar{t}_2y_2, \bar{t}_1z_1, \bar{t}_2z_2; p^{\alpha_1}, p^{\alpha_2})\\
 & \times \widehat{S}_{a, b, d}( \overline{p^{\alpha_1}}x_1, p^{\alpha_2}\overline{p^{\alpha_1}}x_2, p^{\alpha_1}\overline{p^{\alpha_2}}y_1,  \overline{p^{\alpha_2}}y_2, \overline{p^{\alpha_1}}z_1, \overline{p^{\alpha_2}}z_2; t_1 , t_2)  \end{split}
 \end{equation}
 whenever $\alpha_1, \alpha_2 \geq 3$ and $p \nmid t_1t_2$, as well as
 \begin{equation}\label{mult2}
  \begin{split}
 \widehat{S} _{a, b, d}(x_1, x_2, y_1, y_2, z_1, z_2;  t_1u_1, t_2u_2) = &\widehat{S}_{a, b, d}(\bar{t}_1x_1, t_2\bar{t}_1x_2, t_1\bar{t}_2y_1, \bar{t}_2y_2, \bar{t}_1z_1, \bar{t}_2z_2; u_1, u_2)\\
 & \times \widehat{S}_{a, b, d}( \bar{u}_1x_1, u_2\bar{u}_1x_2,  u_1\bar{u}_2y_1, \bar{u}_2y_2, \bar{u}_1z_1,\bar{u}_2z_2; t_1 , t_2)  \end{split}
 \end{equation}
whenever $(t_1t_2, u_1u_2) = 1$. Here $a$ and $b$ on the   right hand sides are understood as primitive residue classes in the respective smaller residue rings.\\ 
 
 
Let $q$ be a prime. By  \eqref{klo6} we have
 \begin{equation}\label{klonew}
\begin{split}
&S^{(N)}(am_1d, bn_2l_2, n_1l_1, m_2d; q^{\alpha_1}, q^{\alpha_2}) \\
= & \sum_{\substack{B_1, C_1 \, ({\rm mod }\, q^{\alpha_1})\\B_2, C_2 \, ({\rm mod }\,  q^{\alpha_2})}} e\left(\frac{am_1dB_1 + n_1l_1(Y_1 q^{\alpha_2}- Z_1 B_2)}{q^{\alpha_1}} + \frac{bn_2l_2B_2 + m_2d(Y_2 q^{\alpha_1} - Z_2B_1)}{q^{\alpha_2}}\right)
\end{split}
\end{equation}
for $N \mid q^{\min(\alpha_1, \alpha_2)}$, 
where the sum is subject to
\begin{equation}\label{cond1}
 q^{\alpha_1}C_2 + B_1B_2 + q^{\alpha_2}C_1 \equiv 0 \, ({\rm mod }\, q^{\alpha_1+\alpha_2}), \quad  (B_j, C_j, q) = 1, \quad  N \mid B_1
\end{equation}
and
\begin{equation}\label{cond2}
   Y_jB_j + Z_jC_j \equiv 1 \, (\text{mod }q^{\alpha_j}) \quad \text{for} \quad  j = 1, 2. 
   \end{equation}
We keep in mind that $q \nmid ab$.  This sum is well-defined as shown in \cite[Lemma 4.2]{BFG}, and does not depend on the choice of the representatives $B_1, B_2$. In particular, we can and will always assume  $$1 \leq B_j \leq q^{\alpha_j} \quad  \text{for}\quad  j = 1, 2.$$  
For future purposes we notice that \eqref{cond1} implies
\begin{equation}\label{lemma1}
  v_q(B_1) \leq \alpha_2, \quad v_q(B_2) \leq \alpha_1
\end{equation}  
where $v_q$ denotes the $q$-adic valuation. Indeed, the first inequality is trivial  if $\alpha_1 \leq \alpha_2$.  If $\alpha_1 > \alpha_2$, the first condition in \eqref{cond1} implies  $q^{\alpha_2} C_1 + B_1 B_2 \equiv 0 \, (\text{mod } q^{\alpha_1})$. If $q \nmid B_1$, there is nothing to prove, otherwise we have $q \nmid C_1$, so $v_q(B_1) \leq v_q(B_1B_2) = \alpha_2$. Similarly one shows the second inequality. 
  
\begin{lemma}\label{peter1} Let  $q$ be a prime, and let $\alpha_1, \alpha_2 \in \Bbb{N}_0$.  For $x_1, y_1, z_1, x_2, y_2, z_2\in \Bbb{Z}$ and $d \in \Bbb{N}$  define $$\gamma := \min\bigl(v_q(x_1), v_q(x_2), v_q(y_1), v_q(y_2), v_q(z_1), v_q(z_2)\bigr), \quad \delta = v_q(d).$$
Then
 we have
$$|\widehat{S}_{a, b, d}(x_1, x_2, y_1, y_2, z_1, z_2; q^{\alpha_1}, q^{\alpha_2})| \leq 3q^{2\min(\alpha_1, \alpha_2)  -  (\alpha_1 + \alpha_2) + 2(\gamma + \delta)}.$$
\end{lemma}  

Here we apply the usual convention $\min(\infty, n) = n$ for $n \in \Bbb{N}_0$. The bound is meaningless for $x_1 = x_2 = y_1 = y_2 = z_1 = z_2 = 0$. This case will be considered in Lemma \ref{zero}. 
We defer the lengthy proof of Lemma \ref{peter1} to the end of this section.   A similar result holds for the twisted transform.

\begin{lemma}\label{twist} Let $\chi$ be a primitive character modulo a prime $p$,   and let $\alpha_1, \alpha_2 \in \Bbb{N}_0$.  For $x_1, x_2, y_1, y_2, z_1, z_2 \in \Bbb{Z}$   define $$\rho := \min\bigl(v_p(x_1), v_p(x_2), v_p(y_1), v_p(y_2), v_p(z_1), v_p(z_2)\bigr).$$  
Assume that $(d, p) = 1$. Then
 we have
$$|\widehat{S}^{\chi}_{a, b, d}(x_1, x_2, y_1, y_2, z_1, z_2; p^{\alpha_1}, p^{\alpha_2})| \leq  3p^{2\min(\alpha_1, \alpha_2)  -  (\alpha_1 + \alpha_2) + 2\rho + 5}.$$
 \end{lemma}  

\textbf{Proof.} 
This follows from the previous lemma  and  the following simple observation. Let $\chi$ be a primitive character modulo $p$, and let $S$ be a $p^{\alpha}$-periodic function with $\alpha \geq 1$. Then
\begin{equation}\label{chartwist}
\frac{1}{p^{\alpha}} \sum_{n \, (\text{mod } p^{\alpha}) } \chi(n) S(n) e\left(-\frac{nx}{p^{\alpha}}\right) = \frac{1}{\tau(\bar{\chi})} \sum_{\beta = 1}^{p-1} \bar{\chi}(\beta) \frac{1}{p^{\alpha}} \sum_{n \, (\text{mod } p^{\alpha}) } S(n) e\left(-\frac{n(x + p^{\alpha-1}\beta)}{p^{\alpha}}\right)
\end{equation}
where as usual $\tau(\chi)$ denotes the Gau{\ss} sum (a complex number of absolute value $p^{1/2}$). We apply this formula   for all six summation  variables in \eqref{shat} and estimate the various $\beta$-sums trivially (this produces an extra factor of $p^{ 6/2} = p^3$). Then we apply Lemma \ref{peter1} with $\gamma \leq \rho + 1$ and $\delta = 0$. 
   
 \begin{lemma}\label{zero} Let $\chi$ be a primitive non-quadratic character modulo a prime $p$, $\alpha_1, \alpha_2 \geq 3$, $(d, p) = 1$. Then 
 $$ \widehat{S}^{\chi}_{a, b, d}(0, 0, 0, 0, 0, 0; p^{\alpha_1}, p^{\alpha_2}) = 0.$$
 \end{lemma} 
 
 \textbf{Proof.} We start with the observation that a Gau{\ss} sum
 $$\sum_{r \, (\text{mod } p^{\lambda})} \chi( r)e\left(\frac{Kr}{p^{\lambda}}\right)$$
 vanishes unless $v_p(K) = \lambda - 1$. 
 
The last condition in \eqref{cond1} implies that $p \mid B_1$, hence $p \nmid C_1$, and we can choose $Y_1 = 0$, $Z_1 = \bar{C}_1$ in \eqref{cond2}.   Hence the $n_1, n_2$-sum becomes
 $$\sum_{n_1 \, (\text{mod }p^{\alpha_1})}\sum_{n_2 \, (\text{mod }p^{\alpha_2})}\bar{\chi}(n_1) \chi(n_2)e\left(-\frac{n_1l_1\bar{C}_1 B_2}{p^{\alpha_1}}\right)e\left(\frac{bn_2l_2B_2}{p^{\alpha_2}}\right).$$
 Since $(l_1l_2, p)=1$ by the presence of the character and $(b, p) = 1$ by assumption, this is only non-zero if $\alpha_1 = \alpha_2 = \alpha \geq 3$, say,  which we assume from now on. Moreover, $v_p(B_2) = \alpha - 1$. Next, the sum over $m_1$ equals
 $$\sum_{m_1 \, (\text{mod } p^{\alpha})} \chi(m_1)e\left(\frac{am_1dB_1}{p^{\alpha}}\right)$$
 which implies $v_p(B_1) = \alpha - 1$. We write $B_1 = p^{\alpha-1} \beta_1$, $B_2 = p^{\alpha -1} \beta_2$ with $p \nmid \beta_1, \beta_2$. Then $C_2 \equiv -C_1 - \beta_1\beta_2 p^{\alpha-2}$ (mod $p^{\alpha}$), and the $B_1, B_2, C_1, C_2$-sum becomes
 \begin{displaymath}
 \begin{split}
&   \sum_{\substack{\beta_1, \beta_2 \, (\text{mod } p)\\ (\beta_1\beta_2, p) = 1}} \sum_{\substack{C_1 \, (\text{mod } p^{\alpha})\\ (C_1, p) = 1}} e\left(\frac{am_1d \beta_1   - n_1l_1\bar{C}_1 \beta_2}{p} + \frac{bn_2l_2\beta_2 + m_2d \overline{(C_1 + \beta_1\beta_2p^{\alpha-2})} \beta_1}{p}\right)\\
= & p^{\alpha-1} \sum_{\substack{\beta_1, \beta_2, C_1 \, (\text{mod } p)\\ (\beta_1\beta_2C_1, p) = 1}}   e\left(\frac{am_1d \beta_1   - n_1l_1\bar{C}_1 \beta_2}{p} + \frac{bn_2l_2\beta_2 + m_2d \bar{C}_1 \beta_1}{p}\right). 
   \end{split}
 \end{displaymath}
 Here we use that $\alpha \geq 3$. 
 The character implies that all variables in the numerator are coprime to $p$. Changing variables $\beta_1 \mapsto \bar{m}_1\beta_1$, $\beta_2\mapsto \overline{n_2l_2}\beta_2$, $C_1\mapsto n_1l_1\overline{n_2l_2}C_1$, we see that this expression depends only on the product $\overline{m_1 n_1l_1}  m_2n_2l_2$, not on the six variables individually. Calling this expression $T(\overline{m_1n_1l_1} m_2n_2l_2)$, we obtain finally by another change of variables (e.g.\ $m_2 \mapsto m_2n_1$) that 
 \begin{displaymath}
 \begin{split}
& \sum_{n_1, n_2, m_1, m_2, l_1, l_2 \, (\text{mod } p^{\alpha})} \bar{\chi}(n_1l_1m_2) \chi(n_2l_2m_1)  T(\overline{m_1n_1l_1 }  m_2n_2l_2)\\
  = &\sum_{n_1, n_2, m_1, m_2, l_1, l_2 \, (\text{mod } p^{\alpha})} \bar{\chi}(n_1^2l_1m_2) \chi(n_2l_2m_1)  T(\overline{m_1l_1}  m_2n_2l_2)
 \end{split}
 \end{displaymath}
and the $n_1$-sum vanishes since $\chi$ is not quadratic.    \\

We combine the previous computations to  the following useful result:
\begin{cor}\label{coro} Let $\chi$ be a primitive non-quadratic character modulo a prime $p$. Let $D_1, D_2, d \in \Bbb{N}$ satisfying $p^3 \mid D_1$, $p^3 \mid D_2$ and $(d, p) = 1$. Let $x_1, x_2, y_1, y_2, z_1, z_2 \in \Bbb{Z}.$ Then
$$|\widehat{S}^{\chi}_{a, b, d}(x_1, x_2, y_1, y_2, z_1, z_2; D_1, D_2)| \leq p^5\tau_3((D_1, D_2)) \frac{(D_1, D_2)^2}{D_1D_2} (x_1, x_2, y_1, y_2, z_1, z_2, D_1, D_2)^2 (d, D_1, D_2)^2$$
where $\tau_3$ is the ternary divisor function. Moreover,  
$$\widehat{S}^{\chi}_{a, b, d}(0, 0, 0, 0, 0, 0; D_1, D_2) =0.$$
\end{cor}

Indeed, the last statement is a direct consequence of   \eqref{mult1} and Lemma \ref{zero}. The first bound follows from the Chinese remainder theorem together with \eqref{mult1}, \eqref{mult2} and Lemmas  \ref{peter1} and \ref{twist}, noting that $\tau_3(q) = 3$ for a prime $q$ and $q^{2\min(\alpha_1, \alpha_2) - (\alpha_1+\alpha_2)} = (q_1^{\alpha_1}, q_2^{\alpha_2})^2/(q_1^{\alpha_1}q_2^{\alpha_2})$. \\

Finally we give the \textbf{proof} of  Lemma \ref{peter1}. We will frequently use the following simple result
\begin{equation}\label{simple}
 \Bigl| \sum_{n, l\, (\text{mod } q^{\alpha})} e\left(\frac{nlB}{q^{\alpha}}\right) e\left(-\frac{nx-lz}{q^{\alpha}}\right) \Bigr|=  \begin{cases} q^{\alpha + v_q(B)}, & v_q(B) \leq \min(v_q(x), v_q(z))\\ 0, & \text{otherwise}\end{cases}
\end{equation}
for integers $x, z$ and $v_q(B) \leq \alpha$.   
Indeed, the $n$-sum vanishes unless $lB - x \equiv 0$ (mod $q^{\alpha}$) in which case it equals $q^{\alpha}$.  This implies in particular $v_q(B) \leq v_q(x)$ and defines $l$ modulo $q^{\alpha - v_q(B)}$. Then the $l$-sum vanishes unless $v_q(z) \geq v_q(B)$, and the result follows easily.

We will now distinguish several cases to estimate $\widehat{S}_{a, b, d}$, but the overall strategy  is always the same. We open the Kloosterman sum and pull the $B_1, B_2, C_1, C_2$-sums outside. Then we sum over $n_1, n_2, l_1, l_2, m_1, m_2$ using orthogonality of additive characters or \eqref{simple}. At this point we estimate trivially and just count how many quadruples $(B_1, B_2, C_1, C_2)$ survive in the outer sum.\\

In order to avoid pathological cases, we treat the case $\alpha_1\alpha_2 = 0$ separately. The case $\alpha_1 = \alpha_2 = 0$ is trivial, so let us assume $\alpha_1 > 0$, $\alpha_2 = 0$. In this case the  $w_6$-Kloosterman sum degenerates to an ordinary Kloosterman sum 
$$S^{(1)} (am_1d, bn_2l_2, n_1l_1, m_2d, q^{\alpha_1}, 1) = S(am_1d, n_1l_1, q^{\alpha_1})$$
by \cite[Property 4.9]{BFG}, so that we need to bound
$$\frac{1}{q^{3\alpha_1}} \sum_{n_1, l_1, m_1 \, (\text{mod } q^{\alpha_1})} S(am_1 d, n_1l_1, q^{\alpha_1}) e\left(\frac{-n_1x_1 - m_1y_1 - l_1z_1}{q^{\alpha_1}}\right).$$
We open the Kloosterman sum, sum over $n_1, l_1$ by \eqref{simple} and over $m_1$ by orthogonality of characters getting the bound $q^{-\alpha_1 + \delta}$. The case $\alpha_1 = 0$, $\alpha_2 > 0$ is similar. From now on we assume $\alpha_1 > 0$ and $\alpha_2 > 0$.\\

We return to the definition \eqref{untwist} and split the Kloosterman sum    in \eqref{klonew} into two parts according to whether $q \nmid C_1$ or $q \mid C_1$. This gives a decomposition  $\widehat{S}_{a, b, d} = T_{a, b, d} + U_{a, b, d}$. The second term where $q \mid C_1$ is easier. In this case $q \nmid B_1$ by \eqref{cond1}, and the first condition there implies   $q \mid B_2$ (since $\alpha_1, \alpha_2 > 0$), so that $q \nmid C_2$. Hence we can choose $$Y_1 = \bar{B}_1,  \quad Z_2 = \bar{C}_2,  \quad Z_1 = Y_2 = 0.$$ 
Considering the $q$-powers in the first condition in \eqref{cond1} again, we see that $v_q(B_1B_2) \leq \alpha_2$, but $v_q(q^{\alpha_2} C_1 + q^{\alpha_1}C_2) \geq \min(\alpha_1, \alpha_2 + 1)$, and we conclude $\alpha_1 \leq \alpha_2$. Thus  
\begin{displaymath}
\begin{split}
U_{a, b, d} = \frac{1}{q^{3(\alpha_1+\alpha_2)}} \sum_{\substack{B_1, C_1 \, (\text{mod }q^{\alpha_1})\\ B_2, C_2 \, (\text{mod }q^{\alpha_2})\\ \eqref{cond1} \text{ holds  with $N=1$}\\ q \nmid B_1C_2, q \mid C_1, q \mid B_2}} & \sum_{\substack{n_1, m_1, l_1 \,(\text{mod }q^{\alpha_1})\\ n_2, m_2, l_2 \,(\text{mod }q^{\alpha_2})}} e\left(\frac{am_1dB_1}{q^{\alpha_1}}+ \frac{bn_2l_2B_2 - m_2d\bar{C}_2 B_1}{q^{\alpha_2}}\right)\\
& \times e\left(-\frac{n_1x_1 + m_1y_1+l_1z_1}{q^{\alpha_1}}\right) e\left(-\frac{n_2x_2 + m_2y_2+l_2z_2}{q^{\alpha_2}}\right). 
\end{split}
\end{displaymath}
We sum trivially over $n_1, l_1$, we use \eqref{simple} in combination with \eqref{lemma1} for the sum over $n_2, l_2$, and we sum over $m_1, m_2$ using orthogonality of characters. The latter two sums leave $q^{\delta}$ choices for $B_1$ and $C_2$ respectively (recall that $q \nmid B_1$).   Now the first condition in \eqref{cond1} determines $B_2$ modulo $q^{\alpha_2}$ which then determines $C_1$.   Altogether we obtain
\begin{equation}\label{u}
|U_{a, b, d}|\leq q^{\alpha_1 - \alpha_2 + 2\delta} = q^{2\min(\alpha_1, \alpha_2) - \alpha_1 - \alpha_2 + 2\delta}.
\end{equation}\\

Now we turn to the estimation of $T_{a, b, d}$ where $q \nmid C_1$, so that 
\begin{equation*}
Z_1 = \bar{C}_1,   \quad Y_1 = 0,
\end{equation*}
and we obtain
\begin{equation}\label{tabd}
\begin{split}
T_{a, b, d} = & \frac{1}{q^{3(\alpha_1+\alpha_2)}} \sum_{\substack{B_1, C_1 \, (\text{mod }q^{\alpha_1})\\ B_2, C_2 \, (\text{mod }q^{\alpha_2})\\ \eqref{cond1} \text{ holds  with $N=1$}\\ q \nmid C_1}}  \sum_{\substack{n_1, m_1, l_1 \,(\text{mod }q^{\alpha_1})\\ n_2, m_2, l_2 \,(\text{mod }q^{\alpha_2})}} e\left(-\frac{n_1x_1 + m_1y_1+l_1z_1}{q^{\alpha_1}}\right) \\
& 
\times e\left(-\frac{n_2x_2 + m_2y_2+l_2z_2}{q^{\alpha_2}}\right) e\left(\frac{am_1dB_1 - n_1l_1\bar{C}_1 B_2}{q^{\alpha_1}}+ \frac{bn_2l_2B_2 + m_2d(Y_2 q^{\alpha_1} - Z_2B_1)}{q^{\alpha_2}}\right).  
\end{split}
\end{equation}
By \eqref{simple}, the sum over $n_2, l_2$ contributes $q^{\alpha_2 + v_q(B_2)}$ if $v_q(B_2) \leq \min(v_q(x_2), v_q(z_2))$ and is zero otherwise. Similarly the sum over $n_1, l_1$ contributes $q^{\alpha_1 + v_q(B_2)}$ if $v_q(B_2) \leq \min(v_q(x_1), v_q(z_1))$ and is zero otherwise (note that by \eqref{lemma1} we have $v_q(B_2) \leq \alpha_1$, so that \eqref{simple} is applicable). We conclude that the combined sum over $n_1, n_2, l_1, l_2$ contributes
\begin{equation}\label{nl}
q^{\alpha_1 + \alpha_2 + 2 \min(v_q(x_1), v_q(x_2), v_q(z_1), v_q(z_2))}.
\end{equation}
  As before   the $m_1$-sum leaves at most $q^{\delta}$ choices for $B_1$, and we fix one of them.  

Let us first assume that $q \nmid B_2$, so that 
\begin{equation}\label{yz1}
Y_2 = \bar{B}_2,  \quad Z_2 = 0.
\end{equation}
  Then the $m_2$-sum leaves at most $q^{\alpha_1 + \delta}$ choices for $B_2$ (and trivially there are at most $q^{\alpha_2})$ choices for $B_2$). Again we fix one of them. If $\alpha_1 \leq \alpha_2$ fix a choice for $C_1$, otherwise fix a choice for $C_2$. In either case, the other $C$-variable is determined by \eqref{cond1}. We conclude that there are in total at most $q^{ \delta + \min(\alpha_1 + \delta, \alpha_2) + \min(\alpha_1, \alpha_2)} \leq q^{2\min(\alpha_1, \alpha_2) + 2\delta}$ choices for the quadruples $(B_1, B_2, C_1, C_2)$ satisfying $q \nmid B_2$.

Let us now assume $q \mid B_2$, so that $q \nmid C_2$ and $$Y_2 = 0,  \quad Z_2 = \bar{C}_2.$$ The $m_2$-sum leaves at most $q^{\delta + v_q(B_1)} \leq q^{\delta + \alpha_1}$ choices for $C_2$ (and trivially there are at most $q^{\alpha_2}$ choices for $C_2$). If $\alpha_1 \leq \alpha_2$ fix a choice for $C_1$, otherwise fix a choice for $B_2$. In either case, the other variable is determined by \eqref{cond1}. As above we conclude that there are in total at most $q^{2\min(\alpha_1, \alpha_2) + 2\delta}$ choices for the quadruples $(B_1, B_2, C_1, C_2)$ satisfying $q \mid B_2$.

We conclude from the previous discussion that the sum over $m_1, m_2$ together with the sum over $B_1, B_2, C_1, C_2$ contributes $2q^{\alpha_1 + \alpha_2 + 2\min(\alpha_1, \alpha_2) + 2\delta}$, and we obtain by \eqref{nl} the total bound
\begin{equation}\label{t1}
|T_{a, b, d}| \leq 2q^{2\min(\alpha_1, \alpha_2)- \alpha_1 - \alpha_2 +  2 \min(v_q(x_1), v_q(x_2), v_q(z_1), v_q(z_2))+ 2\delta}.
\end{equation}\\
 
This is not quite sufficient to substantiate the claim of Lemma \ref{peter1}, so we proceed to prove an alternative bound for $T_{a, b, d}$ as defined in \eqref{tabd}. 
Let us first assume that $q \nmid B_2$, so that \eqref{yz1} holds.  Then by \eqref{simple}, the $n_1, l_1, n_2, l_2$-sums contribute $q^{\alpha_1 + \alpha_2}$, while the $m_1, m_2$-sums leave as above $q^{\delta + \min(\alpha_1 +  \delta, \alpha_2)}$ choices for the pair $(B_1, B_2)$. Fix a choice for $C_1$  if $\alpha_1 \leq \alpha_2$, otherwise fix a choice for $C_2$; in either case the other variable is determined by \eqref{cond1}. In total we obtain at most $q^{2\delta + 2\min(\alpha_1, \alpha_2)}$ choices for the quadruples $(B_1, B_2, C_1, C_2)$, so that together with the $m_1, m_2$-sum we obtain a total contribution of 
\begin{equation}\label{easy}
q^{2 \min(\alpha_1, \alpha_2) - \alpha_1- \alpha_2 + 2\delta}
\end{equation}
for the terms $q \nmid B_2$.

Let us now consider the terms with $q \mid B_2$, so that $q \nmid C_2$ and
$$Y_2 = 0, \quad Z_2 = \bar{C}_2,$$ 
so that \eqref{tabd} simplifies to 
\begin{equation}\label{simplifies}
\begin{split}
\frac{1}{q^{3(\alpha_1+\alpha_2)}} \sum_{\substack{B_1, C_1 \, (\text{mod }q^{\alpha_1})\\ B_2, C_2 \, (\text{mod }q^{\alpha_2})\\ \eqref{cond1} \text{ holds  with $N=1$}\\ q \nmid C_1C_2, q \mid B_2}} & \sum_{\substack{n_1, m_1, l_1 \,(\text{mod }q^{\alpha_1})\\ n_2, m_2, l_2 \,(\text{mod }q^{\alpha_2})}} e\left(\frac{am_1dB_1 - n_1l_1\bar{C}_1 B_2}{q^{\alpha_1}}+ \frac{bn_2l_2B_2 - m_2d\bar{C}_2B_1}{q^{\alpha_2}}\right)\\
& \times e\left(-\frac{n_1x_1 + m_1y_1+l_1z_1}{q^{\alpha_1}}\right) e\left(-\frac{n_2x_2 + m_2y_2+l_2z_2}{q^{\alpha_2}}\right). 
\end{split}
\end{equation}
In the following we assume without loss of generality $1 \leq y_1 \leq q^{\alpha_1}$ and $1 \leq y_2 \leq q^{\alpha_2}$, so that $v_q(y_1) \leq \alpha_1$ and $v_q(y_2) \leq \alpha_2$. It is convenient to first dispense with the case $\min(v_q(y_1), v_q(y_2)) \geq \alpha_2$ (and hence $= \alpha_2$). Here the $n_1, n_2, l_1, l_2$-sums contribute by \eqref{simple} and \eqref{lemma1} at most $q^{\alpha_1 + \alpha_2 + 2\min(\alpha_1, \alpha_2)}$ while there are at most $q^{\delta}$ choices for $B_1$ and trivially at most $q^{2\alpha_2} \leq q^{2\min(v_q(y_1), v_q(y_2))}$ choices for $(B_2, C_2)$ which determines $C_1$. This gives the total bound
\begin{equation}\label{easybound}
q^{2\min(\alpha_1,  \alpha_2) - \alpha_1 - \alpha_2 + 2\min(v_q(y_1), v_q(y_2)) + \delta }
\end{equation}
for \eqref{simplifies} under the present assumption $\min(v_q(y_1), v_q(y_2)) \geq \alpha_2$. From now on we assume 
\begin{equation}\label{eta}
  \min(v_q(y_1), v_q(y_2)) < \alpha_2.
 \end{equation} 
   We distinguish two cases.

Let us first assume $v_q(y_1) \leq v_q(y_2)$. By \eqref{simple} and \eqref{lemma1} the sum over $n_1, l_1, n_2, l_2$ contributes at most $q^{\alpha_1 + \alpha_2 +  2 \min(\alpha_1, \alpha_2)}$. The $m_1$-sum leaves at most $q^{\delta}$ choices for $B_1$ and each of them satisfies $v_q(B_1) \leq v_q(y_1)$. We fix one of them.  Similarly then the sum over $m_2$ leaves at most $q^{\delta + v_q(B_1)} \leq q^{\delta + v_q(y_1)}$ choices for $C_2$. If $B_1, C_2$ are fixed, then the first condition in \eqref{cond1} leaves at most $q^{v_q(B_1)} \leq q^{v_q(y_1)}$  choices for $B_2$ modulo $q^{\alpha_2}$, and the triple $B_1, B_2, C_2$ determines $C_1$. We conclude that there are at most $q^{2v_q(y_1) + 2\delta}$ choices for the quadruple $(B_1, B_2,  C_1, C_2)$, and we obtain the total bound
\begin{equation}\label{case1}
  q^{2 \min(\alpha_1, \alpha_2) - \alpha_1- \alpha_2 + 2 \min(v_q(y_1), v_q(y_2)) + 2\delta} . 
\end{equation}
for \eqref{simplifies} under the present assumption \eqref{eta} and  $v_q(y_1) \leq v_q(y_2)$. This bound dominates \eqref{easybound}. 

The other case $v_q(y_1) > v_q(y_2)$ 
cannot happen: first we observe that the $m_1, m_2$-sum vanishes unless $dB_1 \equiv \bar{a}y_1 \, (\text{mod } q^{\alpha_1})$ and $dB_1 \equiv - C_2y_2 \, (\text{mod } q^{\alpha_2})$. Together with \eqref{eta} this leads to a contradiction unless $v_q(y_2) \geq \alpha_1$. But this is impossible since $ \alpha_1 \geq v_q(y_1) > v_q(y_2)$. 
We summarize that \eqref{case1} is an upper bound for \eqref{simplifies} in all cases, 
and together with \eqref{easy} we conclude
\begin{equation}\label{t2}
|T_{a, b, d}| \leq 2q^{2\min(\alpha_1, \alpha_2)- \alpha_1 - \alpha_2 +  2 \min(v_q(y_1), v_q(y_2) )+ 2\delta}.
\end{equation}
Combining \eqref{t1} and \eqref{t2} with \eqref{u} completes the proof of Lemma \ref{peter1}.

\section{The sixth moment}

\subsection{Setting up the Kuznetsov formula} We prepare now for the proof of Theorem \ref{main}. We recall the setup that $N$ is a large prime, $p$ is a fixed prime and $\chi$ is a primitive non-quadratic character modulo $p$. All implied constants may depend on $p$. Let $\pi \subseteq L_{\text{cusp}}^2(\Gamma_0(N) \backslash \Bbb{H}_3)$ be a cuspidal automorphic representation. Its $L$-function has conductor dividing $N$ \cite[Th\'eor\`eme]{JPSS}. The contribution in the moment estimate of Theorem \ref{main} of those $\pi$ with conductor 1 is independent of $N$, hence $O(1)$, therefore it suffices to consider $\pi$ of conductor $N$. Fix a newvector $\varpi \in \pi$ and denote its normalized Fourier coefficients, defined in \eqref{four} and \eqref{four1}, with $A_{\varpi}(1, 1) = 1$. 
By an approximate functional equation \cite[Theorem 5.2]{IK} we have 
$$|L(1/2, \pi \times \chi)|^2  =  \Bigl|   \sum_n \frac{a_{\pi \times \chi}(n)}{n^{1/2}} V \left(\frac{n}{N^{1/2}}\right) + \eta \sum_n \frac{\overline{a_{\pi \times \chi}(n)}}{n^{1/2}} \overline{V \left(\frac{n}{N^{1/2}}\right)}\Bigr|^2$$
where $a_{\pi \times\chi}(n)$ are the Dirichlet coefficients of $L(s, \pi \times \chi)$, $V$ is a 
 smooth, bounded, rapidly decaying function   depending on $\pi$ and $p$, and $\eta$ is a complex number of absolute value 1 depending on $\pi$ and $p$. The coefficients $a_{\pi\times \chi}(n)$ are multiplicative and satisfy $a_{\pi \times \chi}(n) = A_{\varpi}(n, 1) \chi(n)$ for $(n, Np) = 1$, and $a_{\pi \times \chi}(n) \ll n^{5/14 + \varepsilon}$ by known bounds towards the Ramanujan conjecture on ${\rm GL}(3)$ (although much weaker bounds would suffice for our purpose). Thus we have
 $$L(s, \pi \times \chi) = \sum_{n} \frac{A_{\varpi}(n, 1) \chi(n)}{n^s}  \sum_{\nu=0}^{\infty} \frac{a_{\pi \times \chi}(p^{\nu})}{p^{\nu s}} L_N(s)$$
 for a certain Euler factor  $L_N(s)$.  
 
 We truncate the sums at $n \leq N^{1/2+\varepsilon}$ at the cost of a negligible error. Writing $V$ as its inverse Mellin transform, moving the contour to real part $\varepsilon$ and pulling the rapidly converging integral outside the absolute values, we obtain
$$|L(1/2, \pi \times \chi)|^2  \ll  N^{\varepsilon} \int_{|t| \leq N^{\varepsilon} } \Bigl|\sum_{n\leq N^{1/2+\varepsilon}} \frac{a_{\pi \times \chi}(n)}{n^{1/2+\varepsilon + it}} \Bigr|^2 dt.$$
Coupled with  a smooth partition of unity and the Cauchy-Schwarz inequality we have for some compactly supported  weight functions  $W_j$ (independent of $\pi$) that
\begin{displaymath}
\begin{split}
|L(1/2, \pi \times \chi)|^2 \ll  N^{\varepsilon}  \int_{|t| \leq N^{\varepsilon} }  \Bigl|\sum_{2^j \leq N^{1/2+\varepsilon}} \sum_{n} \frac{a_{\pi \times \chi}(n)}{n^{1/2+it}} W_j\left(\frac{n}{2^j}\right)\Bigr|^2 dt,
\end{split}
\end{displaymath}
up to a negligible error. Since $N$ is prime, the coefficients of the Euler factor $L_N(s)$ are irrelevant  (for $\varepsilon< 1/2$). Estimating the coefficients $a_{\pi \times \chi}(p^{\nu})$ trivially,  we conclude
\begin{displaymath}
\begin{split}
& |L(1/2, \pi \times \chi)|^2   \ll  N^{\varepsilon} \sum_{\nu} \frac{1}{p^{\nu(\frac{1}{2} - \frac{5}{14} -\varepsilon)}} \int_{|t| \leq N^{\varepsilon} }  \Bigl|\sum_{2^j \leq N^{1/2+\varepsilon}} \sum_{n} \frac{A_{\varpi}(n, 1) \chi(n)}{n^{1/2+it}} W_j\left(\frac{np^{\nu}}{2^j}\right)\Bigr|^2 dt
\\&
\ll N^{\varepsilon} \sum_{\nu} \frac{1}{p^{\nu/8}} \int_{|t| \leq N^{\varepsilon} } \sum_{2^j \leq N^{1/2+\varepsilon}} \sum_{n_1, n_2}  \frac{A_{\varpi}(n_1, 1) \overline{A_{\varpi}(n_2, 1)} \chi(n_1)\bar{\chi}(n_2)}{(n_1n_2)^{1/2}}  \left(\frac{n_2}{n_1}\right)^{it}W_j\left(\frac{n_1p^{\nu}}{2^j}\right)\overline{W_j\left(\frac{n_2p^{\nu}}{2^j}\right)}  dt.
\end{split}
\end{displaymath}
We observe that the $n_1, n_2$-sum is non-negative and that (for $\varepsilon< 1/2$) the variables $n_1, n_2$ are coprime to $N$, so that the Fourier coefficients satisfy the unramified Hecke relations, as discussed prior to the statement of Theorem \ref{kuz} (recall that for $\varpi$ with $A_{\varpi}(1, 1) = 0$ all coefficients coming up in the previous sum vanish).   We multiply three such expressions together.  Applying H\"older's inequality to the combined $\nu$-sum and $t$-integral with exponents $2/3$ and $1/3$ and using that
$$\Biggl(\sum_{\nu}  \int_{|t| \leq N^{\varepsilon}} \Bigl(\frac{1}{p^{\nu/8}}\Bigr)^{3/2}\Biggr)^{2/3}   \ll N^{\varepsilon},$$
we obtain
\begin{displaymath}
\begin{split}
  |L(1/2, \pi \times \chi)|^6 \ll N^{\varepsilon} & \sum_{\nu} \frac{1}{p^{\nu/8}}  \int_{|t| \leq N^{\varepsilon} } \sum_{2^j \leq N^{1/2+\varepsilon}}   \sum_{\substack{n_1, m_1, l_1\\ n_2, m_2, l_2}} \chi(n_1m_1l_1)\bar{\chi}(n_2m_2l_2)\left(\frac{n_2m_2l_2}{n_1m_1l_1}\right)^{it}\\
& 
\times \frac{A_{\varpi}(n_1, 1) \overline{A_{\varpi}(n_2, 1)}A_{\varpi}(m_1, 1) \overline{A_{\varpi}(m_2, 1)} A_{\varpi}(l_1, 1) \overline{A_{\varpi}(l_2, 1)}}{(n_1n_2m_1m_2l_1l_2)^{1/2}}  \\
&\times W_j\left(\frac{n_1p^{\nu}}{2^j}\right)\overline{W_j\left(\frac{n_2p^{\nu}}{2^j}\right)} W_j\left(\frac{m_1p^{\nu}}{2^j}\right)\overline{W_j\left(\frac{m_2p^{\nu}}{2^j}\right)} W_j\left(\frac{l_1p^{\nu}}{2^j}\right)\overline{W_j\left(\frac{l_2p^{\nu}}{2^j}\right)}   dt. 
\end{split}
\end{displaymath}
Finally we multiply this with $\sum_{i}| \langle F_i, \tilde{W}_{\mu_{\pi}}\rangle|^2$ where $F_i$ is a collection of functions as in Lemma \ref{omega} and sum over $\pi$. This gives 
\begin{displaymath}
\begin{split}
\sum_{\substack{\pi \subseteq L^2_{\rm cusp}(\Gamma_0(N)\backslash \Bbb{H}_3)\\ \mu_{\pi} \in \Omega}} & |L(1/2, \pi \times \chi)|^6 \ll N^{\varepsilon}  \max_i \max_{|t| \leq N^{\varepsilon} } \max_{M \leq N^{1/2+\varepsilon}} \sum_{ \pi \subseteq L^2_{\rm cusp}(\Gamma_0(N)\backslash \Bbb{H}_3)} | \langle F_i, \tilde{W}_{\mu_{\pi}}\rangle|^2  \\
&\times \Bigl|\sum_{n_1, m_2, l_1} \frac{A_{\varpi}(n_1, 1) \overline{A_{\varpi}(m_2, 1)} A_{\varpi}(l_1, 1) \chi(n_1l_1)\bar{\chi}(m_2)}{(n_1m_2l_1)^{1/2}}   W\left(\frac{n_1}{M}\right) \overline{W\left(\frac{m_2}{M}\right)} W\left(\frac{l_1}{M}\right)\Bigr|^2  
\end{split}
\end{displaymath}
for some smooth compactly supported weight function $W$. 
 In the interest of readable and compact notation  let us introduce  
$$\sum_{m \sim M} f(m) := \sum_{m} f(m) W\left(\frac{m}{M}\right)$$
for some unspecified smooth compactly supported weight function satisfying
$$W^{(j)} \ll_{\varepsilon, j} N^{\varepsilon}$$
for all $j \in \Bbb{N}_0$.   
In other words, $\sim$ has the same meaning as $\asymp$ except that an additional smooth weight function is attached to the sum which comes in handy when one applies Poisson summation. 

We now use the Hecke relation \cite[Theorem 6.4.11]{Go}
$$A_{\varpi}(n, 1) \overline{A_{\varpi}(m, 1)} A_{\varpi}(l, 1)  = \sum_{\substack{d_1d_2 \mid n\\ d_1d_3 \mid m\\ d_2d_3 \mid l}} A_{\varpi}\left(\frac{nl}{d_1d_2^2d_3}, \frac{md_2}{d_1d_3}\right).$$
By \eqref{res} and another application of the Cauchy-Schwarz inequality applied to the $d_1, d_2, d_3$-sum, we obtain 
\begin{displaymath}
\begin{split}
\sum_{\substack{\pi \subseteq L^2_{\rm cusp}(\Gamma_0(N)\backslash \Bbb{H}_3)\\ \mu_{\pi} \in \Omega}}  |L(1/2, \pi \times \chi)|^6&  \ll  \max_i   \max_{M \leq N^{1/2+\varepsilon}}  \frac{N^{2+\varepsilon}}{M^{3}}\sum_{ \pi \subseteq L^2_{\rm cusp}(\Gamma_0(N)\backslash \Bbb{H}_3)} \frac{ | \langle F_i, \tilde{W}_{\mu_{\pi}} \rangle|^2}{\mathcal{N}(\varpi)} \\
&\times \sum_{\substack{d_1, d_2, d_3\\ (d_1d_2d_3, p) = 1}} \Bigl| \sum_{\substack{n \sim M/d_1d_2\\ m \sim M/d_1d_3\\ l \sim M/d_2d_3}}  A_{\varpi}(n l, md_2)  \chi(n l)\bar{\chi}( m )     \Bigr|^2.
\end{split}
\end{displaymath}
By positivity, we can add the rest of the spectrum. For technical reasons it convenient to sum over the spectrum of  $L^2(\Gamma_0(p^3N)\backslash \Bbb{H}_3)$ which contains the sum in the preceding display as oldforms. We open the square and exchange summations. This gives finally our basic inequality 
\begin{equation}\label{basic}
\begin{split}
\sum_{\substack{\pi \subseteq L^2_{\rm cusp}(\Gamma_0(N)\backslash \Bbb{H}_3)\\ \mu_{\pi} \in \Omega}}  |L(1/2, \pi \times \chi)|^6&  \ll  \max_i  \max_{M \leq N^{1/2+\varepsilon}}\frac{N^{2+\varepsilon} }{M^{3}}\sum_{\substack{d_1, d_2, d_3\\(d_1d_2d_3, p) = 1}}  \sum_{\substack{n_1, n_2 \sim M/d_1d_2\\ m_1, m_2 \sim M/d_1d_3\\ l_1, l_2 \sim M/d_2d_3}}   \bar{\chi}(n_1 l_1m_2)\chi( m_1n_2l_2 ) 
 \\
& \times   \int_{(p^3N)} \overline{A_{\varpi}(n_1l_1, m_2d_2)} A_{\varpi}(n_2l_2, m_1d_2)  \frac{ | \langle F_i, \tilde{W}_{\mu_{\pi}} \rangle|^2}{\mathcal{N}(\varpi)}  d\varpi .
\end{split}
\end{equation}

\subsection{Bounding the Kloosterman terms}\label{52}
 The spectral term is now in shape for an application of the Kuznetsov formula (Theorem \ref{kuz}), and accordingly we write the right hand side as a sum of four terms $\Delta + \Sigma_4 + \Sigma_5 + \Sigma_6$, where
 \begin{displaymath}
\begin{split}
 \Delta =  \max_i & \max_{M \leq N^{1/2+\varepsilon}}\frac{N^{2+\varepsilon} }{M^{3}}\sum_{\substack{d_1, d_2, d_3\\(d_1d_2d_3, p) = 1}}  \sum_{\substack{n_1, n_2 \sim M/d_1d_2\\ m_1, m_2 \sim M/d_1d_3\\ l_1, l_2 \sim M/d_2d_3}}   \bar{\chi}(n_1 l_1m_2)\chi( m_1n_2l_2 ) \delta_{\substack{n_1l_1 = n_2l_2\\ m_2= m_1}}  \| F_i\|^2,\\
\Sigma_4 =  \max_i & \max_{M \leq N^{1/2+\varepsilon}}\frac{N^{2+\varepsilon} }{M^{3}}\sum_{\substack{d_1, d_2, d_3\\(d_1d_2d_3, p) = 1}}  \sum_{\substack{n_1, n_2 \sim M/d_1d_2\\ m_1, m_2 \sim M/d_1d_3\\ l_1, l_2 \sim M/d_2d_3}}   \bar{\chi}(n_1 l_1m_2)\chi( m_1n_2l_2 ) \\
&\sum_{\epsilon  = \pm 1} \sum_{\substack{p^3N D_2 \mid D_1\\  d_2m_2 D_1= l_2n_2 D_2^2}}\frac{ \tilde{S}(\epsilon d_2m_1, d_2m_2, l_1n_1, D_2, D_1)}{D_1D_2}  \tilde{\mathcal{J}}_{\epsilon; F_i^{\ast}}\left(\sqrt{\frac{l_1n_1 d_2^2m_1 m_2}{D_1D_2}} \right),\\
\Sigma_5 =\max_i &\max_{M \leq N^{1/2+\varepsilon}}\frac{N^{2+\varepsilon} }{M^{3}}\sum_{\substack{d_1, d_2, d_3\\(d_1d_2d_3, p) = 1}}  \sum_{\substack{n_1, n_2 \sim M/d_1d_2\\ m_1, m_2 \sim M/d_1d_3\\ l_1, l_2 \sim M/d_2d_3}}   \bar{\chi}(n_1 l_1m_2)\chi( m_1n_2l_2 )\\
&\sum_{\epsilon  = \pm 1} \sum_{\substack{ p^3N \mid  D_1 \mid D_2\\ l_1n_1 D_2 = d_2m_1 D_1^2}} \frac{ \tilde{S}(\epsilon l_2n_2, l_1n_1, d_2m_2, D_1, D_2) }{D_1D_2}\tilde{\mathcal{J}}_{\epsilon; F_i}\left( \sqrt{\frac{l_1n_1 l_2n_2 d_2 m_2}{D_1D_2}}  \right), \\
\Sigma_6 =  \max_i & \max_{M \leq N^{1/2+\varepsilon}}\frac{N^{2+\varepsilon} }{M^{3}}\sum_{\substack{d_1, d_2, d_3\\(d_1d_2d_3, p) = 1}}  \sum_{\substack{n_1, n_2 \sim M/d_1d_2\\ m_1, m_2 \sim M/d_1d_3\\ l_1, l_2 \sim M/d_2d_3}}   \bar{\chi}(n_1 l_1m_2)\chi( m_1n_2l_2 )\\
&\hspace{-1cm} \sum_{\epsilon_1, \epsilon_2 = \pm 1} \sum_{p^3N \mid D_1, p^3N \mid D_2  } \frac{S^{(p^3N)}(\epsilon_2 d_2m_1, \epsilon_1 l_2n_2, l_1n_1, d_2m_2, D_1, D_2)}{D_1D_2}\mathcal{J}_{\epsilon;F_i}\left( \frac{\sqrt{l_2n_2 d_2m_2 D_1}}{D_2 }, \frac{\sqrt{l_1n_1 d_2 m_1 D_2}}{D_1}\right).
\end{split}
\end{displaymath}
 
  It is easy to see that $|\Delta| + |\Sigma_4| + |\Sigma_5| \ll N^{2+\varepsilon}$: indeed,   by a divisor argument we have
 $$\Delta \ll  \max_{M \leq N^{1/2+\varepsilon}} \frac{N^{2+\varepsilon} }{M^{3}} \sum_{d_1, d_2, d_3}  \sum_{\substack{n_1, n_2 \sim M/d_1d_2\\ m_1, m_2 \sim M/d_1d_3\\ l_1, l_2 \sim M/d_2d_3}}   \delta_{\substack{n_1l_1 = n_2l_2\\m_2 = m_1}} \ll  \max_{M \leq N^{1/2+\varepsilon}}  \frac{N^{2+\varepsilon}}{M^3} \sum_{d_1, d_2, d_3} \frac{M^3}{(d_1d_2d_3)^2} \ll N^{2+\varepsilon}$$
as desired. For $\Sigma_5$ we simply observe that the conditions $p^3N \mid D_1 \mid D_2$ and $n_1l_1 D_2 = m_1d_2 D_1^2$ implies $N \mid D_1$ and $N^2 \mid D_2$ since $N \nmid n_1l_1$ is prime, so that $D_1D_2 \geq N^3$, but $ n_2l_2m_2d_2n_1l_1 \ll M^{5} \ll N^{5/2+\varepsilon}$, so that (for sufficiently small $\varepsilon > 0$ and sufficiently large $N$) we have $\Sigma_5 = 0$ by Lemma \ref{J1}(a) (with $X_1= X_2 = 1$). Similarly one shows $\Sigma_4 = 0$. \\

 The term corresponding to the long Weyl element can be bounded by
\begin{equation}\label{6}
\begin{split}
& \Sigma_6  \ll \max_{M \leq N^{1/2+\varepsilon}} \frac{N^{1+\varepsilon}}{M^3}  \sum_{\substack{d_1, d_2, d_3\\ (d_2, p) = 1}} \sum_{\epsilon \in \{\pm 1\}^2} 
\Bigl| \sum_{\substack{n_1, n_2 \sim M/d_1d_2\\ m_1, m_2 \sim M/d_1d_3\\ l_1, l_2 \sim M/d_2d_3}}   \bar{\chi}(n_1 l_1m_2)\chi( m_1n_2l_2 ) \\
 &\times \sum_{\substack{p^3\mid D_1\\ p^3 \mid D_2}}  \frac{S^{(p^3)}(\epsilon_2 \bar{N}m_1d_2, \epsilon_1 \bar{N}n_2l_2 , n_1l_1, m_2d_2; D_1, D_2) }{D_1D_2} \mathcal{J}_{\epsilon; F} \left(\frac{\sqrt{n_2l_2m_2d_2D_1}}{ N^{1/2} D_2}, \frac{\sqrt{n_1l_1m_1d_2D_2} }{ N^{1/2} D_1}\right)\Bigr|. 
\end{split}
\end{equation}
where $ \mathcal{J}_{\epsilon; F}$ satisfies the properties of Lemma \ref{J1}(b). 
 Here we used Lemma \ref{prop1}(b) and (c) 
and note that the support of $\mathcal{J}$ given in Lemma \ref{J1}(b) implies
\begin{equation}\label{delta}
D_1, D_2 \ll   \frac{M^3 d_2}{N(d_1d_2d_3)^2} \ll N^{1/2+\varepsilon}.
\end{equation}
so that automatically $(N, D_1 D_2) = 1$. \\

\textbf{Remark:} We pause for a moment and observe that the contribution of the terms $M = N^{1/2}$, $d_1 = d_2 = d_3 = 1$, $D_1 = D_2 \asymp N^{1/2}$ \emph{without} the character $\chi$ exhibits no essential cancellation and is of size $N^{5/2}$ as predicted by the contribution of the maximal Eisenstein series: indeed, the maximal Eisenstein series are parametrized by ${\rm GL}(2)$ cusp forms $f$ for $\Gamma_0(N)$, and a typical Fourier coefficient is given by $A(n, 1) = \sum_{d \mid n} \lambda_f(n)$; thus the maximal Eisenstein contribution is very  roughly of the form
$$\sum_{f} \Bigl| \sum_{n \ll N^{1/2}} \frac{1}{\sqrt{n}} \sum_{d \mid n} \lambda_f(d)\Bigr|^6 \approx \sum_f \Bigl| N^{1/4} L(1, f)\Bigr|^6 \approx N^{5/2}.$$
 In order to obtain the targeted bound $N^2$ we will have to use the extra oscillation of the character.\\

We apply Poisson summation in the 6 variables $n_1, n_2, m_1, m_2, l_1, l_2$.  We call the dual variables $x_1, x_2, y_1, y_2, z_1, z_2$, respectively. By Lemma \ref{J1}(b) the function 
$$J : n_2 \mapsto \mathcal{J}_{\epsilon; F} \left(\frac{\sqrt{n_2l_2m_2d_2D_1}}{ N^{1/2} D_2}, \frac{\sqrt{n_1l_1m_1d_2D_2} }{ N^{1/2} D_1}\right)$$
satisfies
$$n_2^iJ^{(i)}(n_2) \ll_i \left(\frac{M^3d_2}{N (d_1 d_2d_3)^2 D_2}\right)^i$$
for all $i \in \Bbb{N}_0$ under the present size conditions of the variables. We conclude that the dual variable $x_2$ can be bounded by
$$|x_2| \leq N^{\varepsilon} \cdot  \frac{D_2}{M/d_1d_2}\cdot  \frac{M^3d_2}{N (d_1 d_2d_3)^2 D_2} = N^{\varepsilon} \frac{M^2 }{N  d_1d_3^2} \ll N^{\varepsilon}$$
up to a negligible error. By a similar argument, the same bound holds for $x_1$, and we also have
$$|y_1|, |y_2| \leq N^{\varepsilon} \frac{M^2}{N d_1d_2d_3} \ll N^{\varepsilon}, \quad |z_1|, |z_2| \leq N^{\varepsilon} \frac{M^2}{N d_1^2d_3} \ll N^{\varepsilon}.$$
Now we can apply \eqref{new2} with 
$$\alpha_1 = \frac{x_2M}{d_1d_2D_2}, \quad  \beta_1 = \frac{y_2M}{d_1d_3D_2}, \quad \gamma_1 = \frac{z_2M}{d_2d_3D_2}, \quad \alpha_2 = \frac{x_1M}{d_1d_2D_1}, \quad  \beta_2 = \frac{y_1M}{d_1d_3D_1}, \quad \gamma_2 = \frac{z_1M}{d_2d_3D_1}$$
unless   $x_1 x_2y_1y_2z_1z_2 = 0$, in which case we apply \eqref{new1} with $P  =  N^{\varepsilon}$ and 
$$A_1 = \frac{\sqrt{M^3 d_2 D_1}}{d_1d_2d_3N^{1/2}D_2}, \quad A_2 =\frac{\sqrt{M^3 d_2 D_2}}{d_1d_2d_3N^{1/2}D_1} .$$
 In this way we conclude by trivial estimates that 
\begin{displaymath}
\begin{split}
\Sigma_6  \ll \max_{M \leq N^{1/2+\varepsilon}}& \frac{N^{1+\varepsilon}}{M^3}  \sum_{\substack{d_1, d_2, d_3\\ (d_2, p) = 1}} \sum_{\epsilon \in \{\pm 1\}^2} \sum_{p^3 \mid D_1, D_2 \ll \frac{M^3d_2}{N(d_1d_2d_3)^2}} \frac{1}{D_1D_2} \\
& \times  \left(\frac{\min(d_1d_2, d_1d_3, d_2d_3)(D_1D_2)^{1/2}}{M}+\frac{(d_1d_2d_3)^2N(D_1+D_2) }{M^3d_2}\right) \\
&\times \frac{M^6}{(d_1d_2d_3)^4} \sum_{\substack{|x_1|, |x_2| \ll N^{\varepsilon}\\ |y_1|, |y_2| \ll N^{\varepsilon} \\ |z_1|, |z_2| \ll N^{\varepsilon}}} \Bigl|\widehat{S}^{\chi}_{\epsilon_2 \bar{N}, \epsilon_1\bar{N}, d_2}(x_1, x_2, y_1, y_2, z_1, z_2; D_1, D_2)\Bigr|
\end{split}
\end{displaymath}
where $\widehat{S}$ was defined in  \eqref{shat}.  Notice that the $D_1, D_2$-sum restricts the $d_2$-variable to $d_2 \ll N^{1/2+\varepsilon}$. By Corollary \ref{coro} we can bound the innermost sum by $N^{\varepsilon}(D_1, D_2)^{2} d_2^2 (D_1D_2)^{-1}$, and
 obtain
 $$\Sigma_6  \ll  N^{2+\varepsilon} \sum_{\substack{d_1, d_2, d_3\\ d_2 \ll N^{1/2+\varepsilon}}}  \frac{1}{d_1^2d_2d_3^2}\sum_{D_1, D_2 \ll  N^{1/2 + \varepsilon} } \frac{(D_1, D_2)^2(D_1+D_2)}{(D_1D_2)^2 } \ll N^{2+\varepsilon}$$
as desired. This completes the proof of Theorem \ref{main}.  
 

\section{Proofs of Theorems \ref{thm2} - \ref{thm5}} 

For the proof of Theorems \ref{thm2}, \ref{thm3} and \ref{thm5} we choose functions $F_1, \ldots, F_J$ as in Lemma \ref{omega}, and we apply Lemma \ref{J1} with $X_1 = X_2 = 1$. \\ 

For the \textbf{proof of Theorem \ref{thm3}} we proceed as follows. The outer sum is over cuspidal automorphic representations that we interpret as a sum over newvectors. We add artificially the oldforms and the Eisenstein spectrum and bound the mean value in question by
$$ N^{2+\varepsilon}  \sum_{j=1}^J \int_{(N)}  \Bigl|\sum_{\substack{n \asymp X \\ (n, N) = 1}} A_{\varpi}(n, 1) \alpha(n)\Bigr|^2 |\langle F_j, \tilde{W}_{\mu_{\pi} }\rangle|^2  \mathcal{N}(\varpi)^{-1} d\varpi.$$
Here we used also \eqref{res}.  We open the square and apply the Kuznetsov formula. The diagonal term contributes $\ll N^{2+\varepsilon} \| \alpha \|^2$. 
   
 The contribution of the long Weyl element is bounded by
\begin{displaymath}
\begin{split}
&N^{2+\varepsilon} \sum_{\substack{n, m \asymp X\\ (nm, N) = 1}}|\alpha(n)\alpha(m)| \sum_{\substack{ N \mid D_1, D_2\\ D_1, D_2 \ll X}} \frac{|S^{(N)}(\pm 1,  \pm m,  n,  1; D_1, D_2)|}{D_1D_2}   \ll (NX)^{\varepsilon}X^2N^{1/2}\| \alpha \|^2
\end{split}
\end{displaymath}
 by \eqref{average}, since the support condition in Lemma \ref{J1}(b) with $X_1 = X_2 = 1$ restricts $D_1, D_2 \ll X$.

The contribution of the $w_5$ element is bounded by
\begin{displaymath}
\begin{split}
 N^{2+\varepsilon} \sum_{\substack{n, m \asymp X\\ (nm, N) = 1}}|\alpha(n)\alpha(m)| \sum_{\substack{N \mid D_1 \mid D_2\\ n D_2 = D_1^2\\ D_1D_2 \ll X^2}} \frac{|\tilde{S}(\pm m,  n, 1;  D_1, D_2)|}{D_1D_2},\end{split}
\end{displaymath}
again by the support condition in Lemma \ref{J1}(a).  The summation condition implies $D_1 = ndN$, $D_2 = nd^2N^2$, which is only possible if $N \ll 1$, so that we obtain
    \begin{displaymath}
\begin{split}
  \sum_{n, m \asymp X}|\alpha(n)\alpha(m)| \sum_{ d \ll 1}  \frac{|\tilde{S}(\pm m, n, 1; ndN, nd^2N^2)|}{n^2d^3}\ll   \sum_{n, m \asymp X}|\alpha(n)\alpha(m)| \ll X \| \alpha \|^2
  \end{split}
\end{displaymath}
by Lemma \ref{lem1}.   This term is dominated by the $w_6$ contribution. A similar argument works for the $w_4$ contribution. 
This completes the proof.\\
 
 A similar, but simpler, argument shows the bound in Theorem \ref{thm2}.\\
 
 The proof of Theorem \ref{thm5} follows along the lines of \cite[Theorem 2]{BBR}. If $\| \alpha_{\pi}( p)\|_{\infty}    \geq 1 + \delta$, then $\lambda_{\pi}(p^l) \geq (l+1)(l+2)$ for some sufficiently large $l = l(\delta)$, see \cite[(24)]{BBR}. Hence for any $k \geq 1$ we have
 \begin{equation}\label{wehave}
 \begin{split}
   & \#\left\{ \pi \subseteq L_{\text{cusp}}^2(\Gamma_0(N) \backslash \Bbb{H}_3) :  \mu_{\pi} \in \Omega, \, \| \alpha_{\pi}( p) \|  \geq 1 + \delta\right\}  \leq \sum_{\substack{\pi \subseteq L^2_{\text{cusp}}(\Gamma_0(N)\backslash \Bbb{H}_3)\\ \mu_{\pi} \in \Omega}} \frac{|A_{\varpi}(p^l, 1)|^{2k}}{((l+1)(l+2))^{2k}}.
   \end{split}
 \end{equation}
 By \cite[(14)]{BBR} we have
 $$|A_{\varpi}(p^l, 1)|^{2k} = \sum_{r+s \leq 2lk} \alpha_{r, s, l, k} A_{\varpi}(p^r, p^s)A_{\varpi}(1, 1), \quad \sum_{r+s+\leq 2lk} |\alpha_{r, s, l, k}| \leq \left(\frac{(l+1)(l+2)}{2}\right)^{2k}$$
so that by Cauchy-Schwarz and Theorem \ref{thm2} the right hand side of \eqref{wehave} is bounded by
$$\ll (Np^{2lk})^{\varepsilon} 2^{-2k} (N^2 + N^{1/2}p^{2kl})^{1/2}N \ll (Np^{2lk})^{\varepsilon} 2^{-2k} (N^2 + N^{5/4} p^{kl}).$$ Choosing $k = \lfloor \frac{3 \log N}{4l \log p}\rfloor \geq 1$, we obtain
$$   \#\left\{ \pi \subseteq L_{\text{cusp}}^2(\Gamma_0(N) \backslash \Bbb{H}_3) :  \mu_{\pi} \in \Omega, \, \| \alpha_{\pi} (p )\| \geq 1 + \delta\right\}  \ll N^{2 + \varepsilon - \frac{3 \log 2}{2l \log p}}.$$\\

Finally we prove Theorem \ref{thm4}. Here we choose sufficiently large parameters $X_1=  X_2 = X$ to be determined later and apply the Kuznetsov formula with the function $F^{(X, X)}$ as in \eqref{shift}. Then by Lemma \ref{X} and \eqref{res} we have
$$\sum_{\substack{\pi \subseteq L_{\text{cusp}}^2(\Gamma_0(N) \backslash \Bbb{H}_3)\\ \mu_{\pi}  = (\rho + i\gamma, \rho - i\gamma, - 2 i \gamma) \in \Omega\\ |\rho| \geq \varepsilon}} X^{4 +4 |\Re\mu_{\pi}|} \ll N^{2+\varepsilon} \int_{(N)} \frac{|A_{\varpi}(1, 1)|^2 |\langle F^{(X, X)}, \tilde{W}_{\mu_{\pi}}\rangle|^2}{\mathcal{N}(\varpi)} d\varpi.$$
The diagonal term contributes $N^{2+\varepsilon}X^4$. By Lemma \ref{J1}(b), the long Weyl element contributes
$$(NX)^{\varepsilon}N^{2} X^4 \sum_{N \mid D_1, D_2 \ll X^2} \frac{|S^{(N)}(\pm 1, \pm 1,  1,  1; D_1, D_2)|}{D_1D_2}.$$
Assuming $X \leq N^{1-\varepsilon}$ and recalling that $N$ is prime, we have $N^2 \nmid D_1, D_2$. Combining Lemma \ref{prop1}(b), Lemma \ref{prop1}(c) and Lemma \ref{lem1} (with $N = 1$), we obtain the bound
$$(XN)^{\varepsilon} N^{2 } X^4 \sum_{D_1', D_2' \ll  X^2/N} \frac{N(D_1'D_2')^{1/2+\varepsilon}(D_1', D_2')^{1/2}}{D_1'D_2'N^2} \ll (NX)^{\varepsilon} X^6.$$
By Lemma \ref{J1}(a) the $w_5$ element contributes
$$(NX)^{\varepsilon} N^{2 } X^4 \sum_{\substack{N \mid D_1\mid D_2\\ D_2 = D_1^2\\ D_1D_2 \ll X^3}}\frac{\tilde{S}(\pm 1,  1,  1; D_1, D_2)}{D_1D_2} = 0$$
if $X \leq N^{1-\varepsilon}$. Similarly, the $w_4$ contribution vanishes. Choosing $X = N^{1-\varepsilon}$, we obtain the result.

\section*{Acknowledgement} The authors would like to thank the referee for a careful reading of the manuscript.

\end{document}